\RequirePackage{fix-cm}
\documentclass[smallextended]{svjour3}       
\smartqed  
\usepackage{graphicx}%
\usepackage{multirow}%
\usepackage{amsmath,amssymb,amsfonts}%
\usepackage{mathrsfs}%
\usepackage[title]{appendix}%
\usepackage{xcolor}%
\usepackage{textcomp}%
\usepackage{manyfoot}%
\usepackage{booktabs}%
\usepackage{algorithm}%
\usepackage{algorithmicx}%
\usepackage{algpseudocode}%
\usepackage{listings}%
\usepackage{tabularx}
\usepackage{url}

\newtheorem{assumption}{Assumption}%

\newcommand{\R}{\mathbb{R}}
\newcommand{\ip}[2]{\left\langle #1,#2\right\rangle}
\newcommand{\norm}[1]{\left\|#1\right\|}
\newcommand{\D}{\Delta}
\newcommand{\df}{D_f}
\begin{document}

\title{Linear Convergence of Adaptive Accelerated Gradient Methods: A Unified Framework\thanks{This work was partially supported by the Guangdong Major Project of Basic Research [2026B0303000015], the National Natural Science Foundation of China [72495131, 82441027], Guangdong Provincial Key Laboratory of Mathematical Foundations for Artificial Intelligence [2023B1212010001, ZSZDSYS20251011093302004]. The authors declare that they have no competing interests.}
}

\titlerunning{Linear Convergence of Adaptive Accelerated Gradient Methods}        

\author{Weihe Zhang        \and
        Ming Yan 
}


\institute{
           W. Zhang \and M. Yan \at 
           The Chinese University of Hong Kong, Shenzhen\\
           Guangdong, 518172, P.R. China\\
           \email{123090826@link.cuhk.edu.cn, yanming@cuhk.edu.cn}
}

\date{Received: date / Accepted: date}

\maketitle

\begin{abstract}
We establish a unified Lyapunov framework for adaptive accelerated gradient methods on smooth, strongly convex objectives.
Conditions on local curvature, stepsize ratios, and coefficient variation ensure eventual linear convergence without assuming stepsize convergence.
We verify these conditions for specified AdaNAG-G, smooth unconstrained AC-FGM, and AdaAGM policies and use them to design VC-FG and CC-FG.
The verified adaptive main iterations require neither global smoothness nor strong-convexity parameters and avoid restart and backtracking.
The framework also sharpens a fixed-step FISTA/NAG-type bound.
Experiments on four smooth, strongly convex problem families illustrate complementary performance of the new methods.
\end{abstract}

\section{Introduction}
\label{sec:introduction}

We consider the unconstrained minimization problem
\begin{equation}
    \min_{x\in\mathbb{R}^{d}} f(x),
    \label{eq:intro-problem}
\end{equation}
where $f:\mathbb{R}^{d}\to\mathbb{R}$ has an $L$-Lipschitz continuous gradient and is $\mu$-strongly convex, with $0<\mu\le L$.
We denote its unique minimizer by $x^\star$, its optimal value by $f_\star=f(x^\star)$, and its condition number by $\kappa_f=L/\mu$.
This class includes least-squares and logistic regression with positive quadratic regularization.
When the curvature parameters are available, accelerated first-order methods can attain the optimal dependence on $\sqrt{\kappa_f}$ in their iteration complexity~\cite{nesterov2004lectures,daspremont2021acceleration}.
In practice, these parameters may be unavailable or expensive to estimate, and global bounds can be conservative relative to the curvature encountered along an iteration.

Adaptive methods choose stepsizes from function and gradient information gathered along the iteration, without prescribing the global smoothness constant $L$.
They can also use acceleration coefficients that do not require the strong-convexity constant $\mu$.
A sublinear convergence guarantee for general convex objectives, however, does not by itself imply linear convergence under strong convexity.
We seek common conditions that establish eventual linear convergence for different adaptive methods and also guide the construction of new ones.
The analysis allows stepsizes to increase or decrease and does not assume that they converge.

\subsection{Related work}

\paragraph{Curvature adaptation and restarting.}
Backtracking adapts accelerated methods to unknown smoothness, and full-backtracking variants allow stepsizes to increase or decrease~\cite{beck2009fista,scheinberg2014backtracking,calatroni2019backtracking}.
The strongly convex GFISTA analysis uses supplied convexity parameters~\cite{calatroni2019backtracking}.
Universal and adaptive first-order schemes provide further ways to adapt to smoothness~\cite{nesterov2015universal,levy2018accelegrad}.
For unknown strong convexity or growth, adaptive restart schemes exploit the improved geometry~\cite{lin2014adaptive,fercoq2019restart,roulet2020sharpness,aujol2025restart}.
Combining restart with backtracking also permits smoothness adaptation, as in Free-FISTA and RPF-SFISTA~\cite{aujol2024parameterfree,sujanani2025parameterfree}.
Regularization and restarting provide another route to parameter-free guarantees for gradient minimization~\cite{lan2026gradient}.

\paragraph{Linear convergence without restarting.}
Strong convexity can also yield linear convergence for accelerated methods whose updates do not use its strong convexity.
Li, Shi, and Yuan~\cite{li2024linear} establish this property for NAG and FISTA, while Bao et al.~\cite{bao2026restart} prove global linear convergence for original and gradient-restarted accelerated proximal-gradient methods under strong convexity, including composite objectives and more general acceleration sequences.
Geometric convergence under quadratic growth has also been studied with accuracy-dependent inertial parameters~\cite{aujol2024geometrically}.
Related curvature-dependent ISTA and FISTA variants admit linear guarantees when curvature parameters enter their coefficients~\cite{li2026ista}.

Online estimation of strong convexity offers a complementary approach.
The Polyak-type estimate in~\cite{barre2019polyak} uses the known optimal value, while Lyapunov-based schemes update smoothness and convexity estimates with line search when needed~\cite{xu2026adaptive}.
NAG-free~\cite{cavalcanti2025adaptive} estimates strong convexity without restarts, with a global guarantee of gradient-descent order using backtracking; its accelerated guarantee is local and requires a smoothness upper bound and additional assumptions.

\paragraph{Local curvature and acceleration.}
Adaptive gradient and proximal-gradient methods use local curvature estimates and restrictions on stepsize growth to select stepsizes without line search~\cite{malitsky2020adaptive,malitsky2024adaptiveprox}.
AC-FGM~\cite{li2026uniformly} and AdaNAG~\cite{suh2025adanag} combine local curvature adaptation with accelerated convergence guarantees for smooth convex objectives.
The generalized AdaNAG-G family provides further choices of coefficient schedules and stepsize rules.
Borodich and Kovalev~\cite{borodich2026graal} develop accelerated GRAAL, which permits geometric stepsize growth while attaining nearly optimal complexity for smooth convex objectives.
For strongly convex objectives, Wang and Peypouquet~\cite{wang2025adaagm} establish linear convergence of AdaAGM without requiring the strong-convexity modulus or using restart or line search.
Our analysis derives common conditions for linear convergence of specified AdaNAG-G, smooth unconstrained AC-FGM, and AdaAGM policies, without assuming stepsize convergence.

\subsection{Contributions and scope}

\paragraph{A common convergence analysis.}
We bound the change in a common Lyapunov function using a quadratic form in two consecutive gradients.
Positivity conditions for this form, together with suitable lower bounds on the Lyapunov function, lead to contraction.
For stepsizes converging to a positive limit, asymptotic conditions on the coefficients give eventual contraction.
Without assuming stepsize convergence, we obtain iteration-dependent contraction factors.
We then give explicit sufficient conditions on local curvature, ratios of consecutive stepsizes, and coefficient changes under which these factors are eventually bounded above by a common constant smaller than one.
This yields eventual $R$-linear convergence of the objective error and the squared distance to the minimizer.

\paragraph{Verification of existing methods.}
We verify the convergence conditions for specified coefficient and stepsize rules of AdaNAG-G, smooth unconstrained AC-FGM, and AdaAGM.
For AdaAGM, the verification covers the parameter choices for which its published strongly convex theorem gives a positive linear rate~\cite{wang2025adaagm}.
Both the objective error and squared distance to the minimizer admit geometric upper bounds after a finite initial period.
When the acceleration sequence grows polynomially, the bounds also include a factor that decreases as a power of the iteration number.
The contraction rate and the length of the initial period may depend on the objective, initialization, and resulting adaptive trajectory.
The verified adaptive main iterations require neither the strong-convexity modulus nor the global smoothness constant and use no restart.
The precise coefficient rules, parameter restrictions, and initialization conventions are given in Section~\ref{sec:framework-consequences}.

\paragraph{A fixed-step consequence.}
For the smooth FISTA/NAG-type iteration with the prescribed coefficient schedule in Section~\ref{sec:general}, we establish eventual linear convergence for every fixed stepsize $0<s\le1/L$, without using $\mu$ in the updates.
For $0<s<1/L$, setting $a:=\mu s(1-sL)$ gives the bound
\begin{equation*}
    f(x_k)-f_\star=O\!\left(k^{-2}(1+a)^{-k}\right).
    \label{eq:intro-fixed-step-rate}
\end{equation*}
For the same smooth iteration, this improves the factor $(1+a/4)^{-k}$ in Li, Shi, and Yuan~\cite{li2024linear} while retaining the $k^{-2}$ prefactor.
A separate estimate establishes eventual linear convergence at $s=1/L$.

\paragraph{Framework-guided algorithm design.}
We use the convergence conditions to construct variable-coupling and constant-coupling adaptive families.
The families differ in whether an additional weight in the auxiliary gradient update varies with the iteration or remains fixed.
Their stepsize rules bound the product of the stepsize and the local curvature estimate, thereby limiting stepsize growth.
For the selected instances VC-FG and CC-FG, we establish positive lower bounds on the stepsizes and obtain eventual geometric bounds on the objective error and squared distance, with an additional inverse-square factor in the iteration number.
Both methods avoid restart and backtracking.

Section~\ref{sec:setting} presents the unified iteration and the representations of existing methods, and Section~\ref{sec:onestep} derives the common Lyapunov estimates.
Section~\ref{sec:general} treats convergent stepsizes and the fixed-step FISTA/NAG-type specialization.
Section~\ref{sec:nonconvergent} gives contraction conditions without assuming stepsize convergence.
Section~\ref{sec:framework-consequences} verifies these conditions for existing methods and constructs the VC and CC families.
Section~\ref{sec:numerics} reports the experimental comparisons, and Section~\ref{sec:conclusion} concludes.
The appendix provides coefficient bounds used in the verification.

\section{Problem setting and unified iteration}\label{sec:setting}

\subsection{Assumptions and notation}
\label{subsec:setting-notation}

For the minimization problem \eqref{eq:intro-problem}, let $f:\R^d\to\R$ be continuously differentiable. 
The norm $\norm{\cdot}$ and inner product $\ip{\cdot}{\cdot}$ are Euclidean.
Throughout, $f\in\mathcal S^1_{\mu,L}$: there exist $0<\mu\le L$ such that, for all $u,v\in\R^d$,
\begin{equation}\label{eq:strong-convexity}
    \frac\mu2\norm{u-v}^2
    \le f(u)-f(v)-\ip{\nabla f(v)}{u-v}
    \le\frac L2\norm{u-v}^2.
\end{equation}
Thus $f$ is $\mu$-strongly convex and has an $L$-Lipschitz continuous gradient. 
Its unique minimizer is denoted by $x^\star$. We write
\[
    f_\star:=f(x^\star),\qquad
    \Delta_k:=f(x_k)-f_\star,\qquad
    g_k:=\nabla f(x_k).
\]
For $u,v\in\R^d$, the Bregman divergence of $f$ is
\[
    \df(u,v):=f(u)-f(v)-\ip{\nabla f(v)}{u-v}.
\]
We will use the standard bounds
\begin{equation}\label{eq:bregman-gradient-bounds}
    \frac1{2L}\norm{\nabla f(u)-\nabla f(v)}^2
    \le\df(u,v)
    \le\frac1{2\mu}\norm{\nabla f(u)-\nabla f(v)}^2
\end{equation}
and
\begin{equation}\label{eq:PL-and-smooth-gradient}
    \mu^2\norm{x-x^\star}^2
    \le2\mu\bigl(f(x)-f_\star\bigr)
    \le\norm{\nabla f(x)}^2
    \le2L\bigl(f(x)-f_\star\bigr),\qquad x\in\R^d.
\end{equation}

\subsection{The unified three-sequence iteration}
\label{subsec:unified-iteration}

Given a state $(x_k,z_k)$ and coefficients $s_k>0$, $\alpha_k>0$, and $\theta_k,~\theta_{k+1}\ge1$, consider
\begin{equation}\label{eq:scheme}
\begin{cases}
    y_{k+1}=x_k-s_k g_k,\\[2mm]
    z_{k+1}=z_k-s_k\alpha_k\theta_k g_k,\\[2mm]
    x_{k+1}=\left(1-\dfrac1{\theta_{k+1}}\right)y_{k+1}
             +\dfrac1{\theta_{k+1}}z_{k+1}.
\end{cases}
\end{equation}
The last line of \eqref{eq:scheme} gives
\begin{equation}\label{eq:z-identity}
    z_{k+1}=x_{k+1}+(\theta_{k+1}-1)(x_{k+1}-y_{k+1}).
\end{equation} 
Method-specific initialization and coefficient policies are given in Section~\ref{sec:framework-consequences}. 

\subsection{Local curvature estimate}
\label{subsec:local-curvature}

Once $x_{k+1}$ and $g_{k+1}$ are available, define
\begin{equation}\label{eq:Lk-bregman-form}
    L_{k+1}:=
    \begin{cases}
    \dfrac{\norm{g_{k+1}-g_k}^2}{2\df(x_k,x_{k+1})},
        &x_{k+1}\ne x_k,\\[4mm]
    L_k,&x_{k+1}=x_k.
    \end{cases}
\end{equation}
This quantity can then be used to select the next stepsize. 
Strong convexity makes the denominator positive for distinct iterates; the second case specifies a positive continuation for coincident iterates.

A positive initial estimate in $[\mu,L]$ can be obtained without knowing either endpoint. 
For any two distinct points $\bar x_0,\widetilde x_0$, set
\[
    L_0:=\frac{\norm{\nabla f(\bar x_0)-\nabla f(\widetilde x_0)}}
                    {\norm{\bar x_0-\widetilde x_0}}.
\]
Strong monotonicity and Lipschitz continuity of $\nabla f$ give $\mu\le L_0\le L$. 
A method may use a separate secant estimate during initialization, as AC-FGM does. 

\begin{lemma}[Bounds and identity for the local curvature]\label{lem:Lk-bounds}
Suppose that the preceding curvature estimate lies in $[\mu,L]$.
For every $k\ge0$ at which $L_{k+1}$ is defined by
\eqref{eq:Lk-bregman-form},
\begin{equation}\label{eq:Lk-interval}
    \mu\le L_{k+1}\le L
\end{equation}
and
\begin{equation}\label{eq:Lk-bregman-identity}
    \df(x_k,x_{k+1})
    =\frac1{2L_{k+1}}\norm{g_{k+1}-g_k}^2.
\end{equation}
\end{lemma}

\begin{proof}
If $x_{k+1}\ne x_k$, the two bounds in~\eqref{eq:bregman-gradient-bounds} give \eqref{eq:Lk-interval}, and~\eqref{eq:Lk-bregman-identity} follows from the definition. 
If $x_{k+1}=x_k$, both sides of equation~\eqref{eq:Lk-bregman-identity} vanish and the interval bound is inherited from $L_k$.
\end{proof}

\subsection{Representative methods}
\label{subsec:represented-methods}

\subsubsection{AdaNAG-G}\label{adanagg-scheme}

AdaNAG-G~\cite{suh2025adanag} already has the form \eqref{eq:scheme}.
Its prescribed coefficient sequences and adaptive stepsize rule~\eqref{eq:adanagg-stepsize} are specified, along with their verification, in Section~\ref{sec:framework-consequences}.

\subsubsection{AC-FGM}\label{sec:acfgm-setting}

In the smooth unconstrained case, the AC-FGM updates of Li and Lan~\cite{li2026uniformly} can be written as
\begin{equation*}
\begin{cases}
    y_{k+1}^{\rm AC}=z_k-\eta_{k+1}g_k,\\[2mm]
    z_{k+1}=z_k-\beta^{\rm AC}\eta_{k+1}g_k,\\[2mm]
    x_{k+1}=\dfrac{\tau_{k+1}^{\rm AC}x_k+y_{k+1}^{\rm AC}}
                      {1+\tau_{k+1}^{\rm AC}},
\end{cases}
\end{equation*}
where $\eta_{k+1}>0$ and $0<\beta^{\rm AC}<1$.
For $\tau_k^{\rm AC}>0$ and $\tau_{k+1}^{\rm AC}>0$, set
\begin{equation}\label{eq:acfgm-unified-parameters}
    \theta_k=1+\tau_k^{\rm AC},\qquad
    \alpha_k=\frac{\beta^{\rm AC}\tau_{k+1}^{\rm AC}}
                       {(1-\beta^{\rm AC})(1+\tau_k^{\rm AC})},\qquad
    s_k=\frac{(1-\beta^{\rm AC})\eta_{k+1}}{\tau_{k+1}^{\rm AC}}.
\end{equation}
These choices give
\(s_k\alpha_k\theta_k=\beta^{\rm AC}\eta_{k+1},~
\tau_{k+1}^{\rm AC}s_k=(1-\beta^{\rm AC})\eta_{k+1}\).
Consequently, the same $x_k,z_k$ satisfy \eqref{eq:scheme} with the unified intermediate point $y_{k+1}=x_k-s_kg_k$. 
This point generally differs from $y_{k+1}^{\rm AC}$. 
The policy used in Section~\ref{sec:framework-consequences} enters the unified recursion from $k=2$ after a separate initialization.

\subsubsection{AdaAGM}\label{sec:adaagm-setting}

The AdaAGM representation of Wang and Peypouquet~\cite{wang2025adaagm} has the form
\begin{equation}\label{eq:adaagm-update}
\begin{cases}
    y_{k+1}=x_k-s_kg_k,\\[2mm]
    x_{k+1}=y_{k+1}
      +\dfrac{\theta_k-1}{\theta_{k+1}}(y_{k+1}-y_k)
      +\dfrac{(\gamma-1)\theta_k}{\theta_{k+1}}(y_{k+1}-x_k),\\[2mm]
    s_{k+1}=T_k(s_k,L_{k+1}),
\end{cases}
\end{equation}
where $\gamma>0$ and $T_k$ is a chosen stepsize policy. 
Define
\[
    z_k:=\theta_kx_k-(\theta_k-1)y_k.
\]
The second update in \eqref{eq:adaagm-update} then yields $z_{k+1}=z_k-\gamma s_k\theta_kg_k$. 
Thus AdaAGM is represented by~\eqref{eq:scheme} with $\alpha_k\equiv\gamma$.
Section~\ref{sec:framework-consequences} specifies the coefficient sequences, adaptive stepsize rule, and parameter restrictions used in our verification.

\section{Common Lyapunov estimates}\label{sec:onestep}

The shared estimates combine strong convexity and the local curvature identity. 
We leave the energy weights free for the convergence analyses in Sections~\ref{sec:general} and~\ref{sec:nonconvergent}.
For $k\ge1$, consider
\begin{equation}\label{eq:general-energy}
    \mathcal E_k:=\norm{z_k-x^\star}^2
       +B_k\norm{g_{k-1}}^2+C_k\Delta_{k-1},
\end{equation}
where $\{B_k\}$ and $\{C_k\}$ are two real sequences.
This template permits $B_k<0$. 

Choose an auxiliary sequence $\beta_k\in(0,1)$ and define
\begin{equation}\label{eq:common-quadratic-F}
\begin{aligned}
    W_k:=&\alpha_k\theta_k(\theta_k-1),\quad
    H_k:=\alpha_k^2\theta_k^2,\\
    F_k:=&2s_kW_k\bigl[1-\mu s_{k-1}(1-\beta_k)\bigr].
\end{aligned}
\end{equation}
The weights $B_k,C_k$ and the balancing parameter $\beta_k$ may involve $\mu$ or $L$. 
In particular, $\beta_k$ is distinct from the algorithmic parameters $\beta^{\rm AC}$ and $\beta^{\rm AGM}$.
Throughout this section, $\theta_k>1$ at the indices under consideration.

\subsection{The common one-step estimate}
\label{subsec:common-one-step}
\begin{lemma}[Strong-convexity cross-step estimate]\label{lem:cross-step}
Let $f\in\mathcal S^1_{\mu,L}$. At any index $k\ge1$ for which $y_k=x_{k-1}-s_{k-1}g_{k-1}$ with $s_{k-1}>0$,
\begin{equation}\label{eq:cross-step}
\begin{aligned}
    &(1-\mu s_{k-1})\bigl(f(x_k)-f(x_{k-1})\bigr)
       -\ip{g_k}{x_k-x_{k-1}}\\
    &\qquad\le\frac{\mu s_{k-1}^2}{2}\norm{g_{k-1}}^2
                -\frac\mu2\norm{x_k-y_k}^2.
\end{aligned}
\end{equation}
\end{lemma}

\begin{proof}
Strong convexity at $x_k$ gives
\[
    f(x_{k-1})\ge f(x_k)+\ip{g_k}{x_{k-1}-x_k}
                         +\frac\mu2\norm{x_{k-1}-x_k}^2.
\]
Using $x_{k-1}=y_k+s_{k-1}g_{k-1}$, we obtain
\begin{align*}
    &f(x_k)-f(x_{k-1})-\ip{g_k}{x_k-x_{k-1}}\\
    &\quad\le-\frac\mu2\norm{x_k-x_{k-1}}^2\\
    &\quad=-\frac\mu2\norm{x_k-y_k}^2
          +\frac{\mu s_{k-1}^2}{2}\norm{g_{k-1}}^2
          +\mu s_{k-1}\ip{g_{k-1}}{x_k-x_{k-1}}.
\end{align*}
Convexity of $f$ implies $\ip{g_{k-1}}{x_k-x_{k-1}}\le f(x_k)-f(x_{k-1})$.
Multiplying this inequality by $\mu s_{k-1}>0$ and rearranging proves~\eqref{eq:cross-step}. 
No restriction on the sign of $1-\mu s_{k-1}$ is needed.
\end{proof}

\begin{lemma}[One-step Lyapunov estimate]\label{lem:onestep}
Let $f\in\mathcal S^1_{\mu,L}$. At an index $k\ge1$, suppose that
the previous gradient step and coupling identity hold:
\[
    y_k=x_{k-1}-s_{k-1}g_{k-1},\qquad
    z_k=x_k+(\theta_k-1)(x_k-y_k),\qquad s_{k-1}>0.
\]
Suppose also that the current update~\eqref{eq:scheme} and the
curvature identity~\eqref{eq:Lk-bregman-identity} for
$(x_{k-1},x_k)$ hold. Then the energy \eqref{eq:general-energy} satisfies
\begin{align}\label{eq:onestep}
\mathcal E_{k+1}-\mathcal E_k\le{}&
    (C_{k+1}-F_k-2s_k\alpha_k\theta_k)\Delta_k
       +(F_k-C_k)\Delta_{k-1}\nonumber\\
    &-\mu s_k\alpha_k\theta_k\norm{x_k-x^\star}^2
       -\mu s_kW_k(1-\beta_k)\norm{x_k-y_k}^2\nonumber\\
    &+(B_{k+1}+s_k^2H_k)\norm{g_k}^2+\bigl[-B_k+\mu s_kW_k(1-\beta_k)s_{k-1}^2\bigr]
          \norm{g_{k-1}}^2\nonumber\\
    &-2s_{k-1}s_kW_k\ip{g_k}{g_{k-1}}
       -\frac{\beta_ks_kW_k}{L_k}\norm{g_k-g_{k-1}}^2.
\end{align}
The two negative distance terms can be replaced by
\begin{equation}\label{eq:onestep-z-distance}
    -\frac{\mu s_k\alpha_k\theta_k(1-\beta_k)}{\theta_k-\beta_k}
          \norm{z_k-x^\star}^2.
\end{equation}
\end{lemma}

\begin{proof}
The $z$-update and \eqref{eq:z-identity} give
\begin{align}
    &\norm{z_{k+1}-x^\star}^2-\norm{z_k-x^\star}^2
    =-2s_k\alpha_k\theta_k\ip{g_k}{z_k-x^\star}
      +s_k^2H_k\norm{g_k}^2 \nonumber\\
    &\qquad=-2s_k\alpha_k\theta_k\ip{g_k}{x_k-x^\star}
      -2s_kW_k\ip{g_k}{x_k-y_k}
      +s_k^2H_k\norm{g_k}^2. \label{eq:z-square-expanded}
\end{align}
Strong convexity yields
\begin{equation}\label{eq:first-inner-bound}
    \ip{g_k}{x_k-x^\star} \ge \Delta_k+\frac\mu2\norm{x_k-x^\star}^2.
\end{equation}
To estimate the other inner product, write $x_k-y_k=x_k-x_{k-1}+s_{k-1}g_{k-1}$. 
The Bregman curvature identity~\eqref{eq:Lk-bregman-identity} gives the exact relation
\begin{equation*}
    \ip{g_k}{x_k-x_{k-1}} =\Delta_k-\Delta_{k-1}
    +\frac1{2L_k}\norm{g_k-g_{k-1}}^2.
\end{equation*}
Combine this identity with the lower bound on the same inner product from Lemma~\ref{lem:cross-step}, using weights $\beta_k$ and $1-\beta_k$, respectively. 
Adding $s_{k-1}\ip{g_k}{g_{k-1}}$ gives
\begin{align}\label{eq:second-inner-bound}
\ip{g_k}{x_k-y_k} \ge{}& \bigl[1-\mu s_{k-1}(1-\beta_k)\bigr] (\Delta_k-\Delta_{k-1}) \nonumber\\
&+\frac{(1-\beta_k)\mu}{2}\norm{x_k-y_k}^2
 -\frac{(1-\beta_k)\mu s_{k-1}^2}{2}\norm{g_{k-1}}^2\nonumber\\
&+s_{k-1}\ip{g_k}{g_{k-1}}
 +\frac{\beta_k}{2L_k}\norm{g_k-g_{k-1}}^2.
\end{align}
Substitute \eqref{eq:first-inner-bound} and~\eqref{eq:second-inner-bound} into \eqref{eq:z-square-expanded}, and add \(B_{k+1}\norm{g_k}^2-B_k\norm{g_{k-1}}^2 +C_{k+1}\Delta_k-C_k\Delta_{k-1}\). 
The definition of $F_k$ then proves \eqref{eq:onestep}.

For the new distance term in~\eqref{eq:onestep-z-distance}, use
\[
    a\|u\|^2+b\|v\|^2\ge\frac{ab}{a+b}\|u+v\|^2,
    \qquad a,b>0,
\]
with $u=x_k-x^\star$, $v=(\theta_k-1)(x_k-y_k)$, $a=s_k\alpha_k\theta_k$, and $b=s_k\alpha_k\theta_k(1-\beta_k)/(\theta_k-1)$. 
Since $u+v=z_k-x^\star$, we obtain
\begin{align*}
&s_k\alpha_k\theta_k\norm{x_k-x^\star}^2
 +s_kW_k(1-\beta_k)\norm{x_k-y_k}^2\nonumber\\
 \ge&
\frac{s_k\alpha_k\theta_k(1-\beta_k)}{\theta_k-\beta_k}\norm{z_k-x^\star}^2.
\end{align*}
Multiplying by $\mu$ gives~\eqref{eq:onestep-z-distance}.
\end{proof}

\subsection{Reduction to a quadratic form}
\label{subsec:common-quadratic-reduction}

Choosing $C_k=F_k$ cancels the previous objective-gap term in~\eqref{eq:onestep}. 
Any positive coefficient left on $\Delta_k$ can be absorbed into $\norm{g_k}^2$ using~\eqref{eq:PL-and-smooth-gradient}. 
The remaining interaction between $g_k$ and $g_{k-1}$ is a quadratic form whose coefficients depend on the stepsizes, acceleration coefficients, and energy weights.

\begin{lemma}[Common quadratic reduction]\label{lem:common-quadratic}
Under the hypotheses of Lemma~\ref{lem:onestep}, take $C_k=F_k$ in
\eqref{eq:general-energy}. Define
\begin{equation}\label{eq:common-quadratic-residual}
    R_k:=\frac{[F_{k+1}-F_k-2s_k\alpha_k\theta_k]_+}
    {2\mu s_k^2W_k},
    \qquad [a]_+:=\max\{a,0\},
\end{equation}
and the coefficients
\begin{equation}\label{eq:common-quadratic-coefficients}
\begin{aligned}
    A_k&:=\frac{\beta_k}{s_kL_k}
       -\frac{B_{k+1}}{s_k^2W_k}-\frac{H_k}{W_k}-R_k,\\
    D_k&:=\frac{s_{k-1}}{s_k}-\frac{\beta_k}{s_kL_k},\\
    G_k&:=\frac{\beta_k}{s_kL_k}
       +\frac{B_k}{s_k^2W_k}
       -\mu(1-\beta_k)\frac{s_{k-1}^2}{s_k}.
\end{aligned}
\end{equation}
Set
\begin{equation*}
    Q_k:=A_k\norm{g_k}^2
       +2D_k\ip{g_k}{g_{k-1}}+G_k\norm{g_{k-1}}^2.
\end{equation*}
Then, at every such index $k\ge1$,
\begin{equation}\label{eq:common-quadratic-z-descent}
    \mathcal E_{k+1}-\mathcal E_k
    \le-\frac{\mu s_k\alpha_k\theta_k(1-\beta_k)}{\theta_k-\beta_k}\norm{z_k-x^\star}^2-s_k^2W_kQ_k.
\end{equation}
\end{lemma}

\begin{proof}
With $C_k=F_k$, the coefficient of $\Delta_{k-1}$ in~\eqref{eq:onestep} vanishes. Since $\Delta_k\ge0$ and $2\mu\Delta_k\le\norm{g_k}^2$,
\begin{equation*}
    (F_{k+1}-F_k-2s_k\alpha_k\theta_k)\Delta_k
    \le R_ks_k^2W_k\norm{g_k}^2.
\end{equation*}
Expanding $\norm{g_k-g_{k-1}}^2$ and collecting the three gradient terms in~\eqref{eq:onestep} gives \eqref{eq:common-quadratic-z-descent} with the coefficients in \eqref{eq:common-quadratic-coefficients}.
\end{proof}

The role of the quadratic form is especially clear when $A_k>0$:
\[
    Q_k=A_k\norm{g_k+\frac{D_k}{A_k}g_{k-1}}^2
           +\left(G_k-\frac{D_k^2}{A_k}\right)\norm{g_{k-1}}^2.
\]
A positive Schur complement $G_k-D_k^2/A_k$ therefore yields an energy decrease proportional to the squared norm of the previous gradient. 
Combining this decrease with the negative distance term and suitable bounds on the energy gives the contraction estimates in the next two sections.

\section{The convergent-stepsize regime}\label{sec:general}
We consider stepsizes $s_k\to s>0$ and use the notation $W_k,H_k$ from \eqref{eq:common-quadratic-F}. 
We impose the following asymptotic conditions on the coefficients and stepsizes. 
These conditions require neither convergence of the local curvature estimates $L_k$ nor monotonicity of the stepsizes $s_k$.

\begin{assumption}[Asymptotic conditions]\label{ass:general-structure}
The parameters in \eqref{eq:scheme} satisfy
\begin{enumerate}
    \item $\liminf\limits_{k\to\infty}\alpha_k=\underline\alpha>0$,
    $\theta_k\to\infty$, $W_{k+1}/W_k\to1$, and
    $\limsup\limits_{k\to\infty}H_k/W_k\le1$.
    \item $s_k\to s>0$ and $\limsup\limits_{k\to\infty}s_kL_k\le1$.
\end{enumerate}
\end{assumption}

\begin{lemma}[Vanishing normalized residual]\label{lem:residual-absorption}
Fix $\beta_k\equiv\beta\in(0,1)$ in Lemma~\ref{lem:common-quadratic}. 
If $s_k\to s>0$, $\theta_k\to\infty$, and $W_{k+1}/W_k\to1$, then $R_k\to0$.
\end{lemma}

\begin{proof}
The definition of $F_k$ gives
$F_k/W_k\to2s[1-\mu s(1-\beta)]$.
Together with $W_{k+1}/W_k\to1$, this implies
\[
    \frac{F_{k+1}-F_k}{s_k^2W_k}\longrightarrow0, \qquad
    \frac{2s_k\alpha_k\theta_k}{s_k^2W_k}
       =\frac{2}{s_k(\theta_k-1)}\longrightarrow0.
\]
Taking the positive part in~\eqref{eq:common-quadratic-residual} yields $R_k\to0$.
\end{proof}

\subsection{Eventual contraction for convergent stepsizes}

\begin{theorem}[Eventual linear convergence for convergent stepsizes]\label{thm:general-linear}
Let $f\in\mathcal S^1_{\mu,L}$. 
Suppose that, for some integer $k_0\ge1$, the update and curvature identities required by Lemma~\ref{lem:onestep} hold for every $k\ge k_0$. 
Assume also that Assumption~\ref{ass:general-structure} holds.
Write
\[
    \bar q:=\limsup_{k\to\infty}s_kL_k,\qquad m:=\mu s.
\]
Then $0<m\le\bar q\le1$ and $0<\underline\alpha\le1$.

\smallskip
\noindent\textnormal{(i) The case $\bar q<1$.}
Fix $\beta\in(\bar q,1)$ and define
\begin{equation*}
    V_k^{\mathrm{loc},\beta}:=\norm{z_k-x^\star}^2
       +2s_kW_k[1-\mu s_{k-1}(1-\beta)]\Delta_{k-1}.
\end{equation*}
For every $0<\eta<m\underline\alpha(1-\beta)$, there is $N\ge k_0$
such that, for all $k\ge N$,
\begin{align}
    V_{k+1}^{\mathrm{loc},\beta}
       &\le(1-\eta)V_k^{\mathrm{loc},\beta},
          \label{eq:thm42-local-contraction}\\
    V_k^{\mathrm{loc},\beta}
       &\ge\norm{z_k-x^\star}^2+2\beta s_kW_k\Delta_{k-1}.
          \label{eq:thm42-local-lower}
\end{align}

\smallskip
\noindent\textnormal{(ii) The case $sL>1/2$.}
Set
\begin{equation*}
    p:=\frac1{2Ls}\in(0,1),\qquad
    \varepsilon:=\frac{p(1-p)}{2+pm},\qquad
    \beta:=1-\varepsilon\in(7/8,1),
\end{equation*}
and define
\begin{equation*}
\begin{aligned}
    V_k^{\mathrm{glob}}:= {}&\norm{z_k-x^\star}^2
       -\frac{s_kW_k}{2L}\norm{g_{k-1}}^2\\
       &+2s_kW_k(1-\mu s_{k-1}\varepsilon)\Delta_{k-1}.
\end{aligned}
\end{equation*}
For every $0<\eta<m\underline\alpha\varepsilon$, there is $N\ge k_0$ such that, for all $k\ge N$,
\begin{align}
    V_{k+1}^{\mathrm{glob}}
       &\le(1-\eta)V_k^{\mathrm{glob}},
          \label{eq:thm42-global-contraction}\\
    V_k^{\mathrm{glob}}
       &\ge\norm{z_k-x^\star}^2+\frac34s_kW_k\Delta_{k-1}.
          \label{eq:thm42-global-lower}
\end{align}
\end{theorem}

\begin{proof}
Since $L_k\ge\mu$ and $s_k\to s$, we have $0<m\le\bar q\le1$. Also,
\[
    \frac{H_k}{W_k}=\frac{\alpha_k\theta_k}{\theta_k-1},
\]
so (S1) implies $\limsup\limits_{k\to\infty}\alpha_k\le1$ and hence $0<\underline\alpha\le1$. 
For either part, take $\beta_k\equiv\beta$. Lemma~\ref{lem:residual-absorption} gives $R_k\to0$, and the coefficient multiplying $-\norm{z_k-x^\star}^2$ in~\eqref{eq:common-quadratic-z-descent} satisfies
\begin{align}\label{eq:xi_k_lower}
    \liminf_{k\to\infty}\frac{\mu s_k\alpha_k\theta_k(1-\beta)}{\theta_k-\beta}
    \ge m\underline\alpha(1-\beta).
\end{align}
For each fixed $\beta$, put $t_k:=\beta/(s_kL_k)$. 
The sequence $\{t_k\}$ is bounded because $s_k\to s>0$ and $L_k\ge\mu$.
We also have $s_{k-1}/s_k\to1$ and $H_k/W_k\le1+o(1)$.

\smallskip
\noindent\emph{Part (i).}
Apply Lemma~\ref{lem:common-quadratic} with $B_k=0$ and $\beta_k\equiv\beta$. 
Since \(\liminf\limits_{k\to\infty}\beta/(s_kL_k)= \beta/\bar q>1\), we have \(\liminf\limits_{k\to\infty}A_k\ge \beta/\bar q-1>0\). 
The coefficient formulas give
\[
    A_k\ge t_k-1-o(1),\qquad
    D_k=1-t_k+o(1),\qquad
    G_k=t_k-m(1-\beta)+o(1).
\]
Since $t_k$ is bounded and $t_k-1$ is bounded away from zero,
completing the square and using
$t-(1-t)^2/(t-1)=1$ yields
\begin{align*}
    \liminf_{k\to\infty}\left(G_k-\frac{D_k^2}{A_k}\right)
    &\ge1-m(1-\beta)=:b>0.
\end{align*}
For every $\gamma\in(0,b)$ and all sufficiently large $k$,
\begin{equation}\label{eq:thm42-local-descent}
    V_{k+1}^{\mathrm{loc},\beta}-V_k^{\mathrm{loc},\beta}
    \le-\frac{\mu s_k\alpha_k\theta_k(1-\beta)}{\theta_k-\beta}\|z_k-x^\star\|^2
       -\gamma s_k^2W_k\|g_{k-1}\|^2.
\end{equation}

Fix $0<\eta<m\underline\alpha(1-\beta)$. Since $\eta<m$, choose $\gamma= (\eta +m)b/(2m)$. Then $\eta b/m<\gamma<b$, and
\[
    \frac{\gamma\mu s_k}{1-\mu s_{k-1}(1-\beta)} \to\frac{\gamma m}{b}>\eta.
\]
Using $\norm{g_{k-1}}^2\ge2\mu\Delta_{k-1}$, we obtain eventually
\begin{equation*}
\begin{aligned}
    \gamma s_k^2W_k\|g_{k-1}\|^2
    &\ge2\gamma\mu s_k^2W_k\Delta_{k-1}\\
    &\ge2\eta s_kW_k[1-\mu s_{k-1}(1-\beta)]\Delta_{k-1}.
\end{aligned}
\end{equation*}
Together with \eqref{eq:xi_k_lower}, this bounds the right-hand side of \eqref{eq:thm42-local-descent} by $-\eta V_k^{\mathrm{loc},\beta}$ and proves \eqref{eq:thm42-local-contraction}.
Finally, $m\le\bar q<1$ implies $\mu s_{k-1}<1$ eventually, so $1-\mu s_{k-1}(1-\beta)\ge\beta$. This proves~\eqref{eq:thm42-local-lower}.

\smallskip
\noindent\emph{Part (ii).}
Take $B_k=-s_kW_k/(2L)$ and $\beta_k\equiv1-\varepsilon$.
The assumptions give
\[
    \frac{-B_{k+1}}{s_k^2W_k}\longrightarrow p,
    \qquad \frac{-B_k}{s_k^2W_k}\longrightarrow p.
\]
Also $0<\varepsilon<p$ and $\varepsilon<1/8$. 
The coefficient formulas give
\[
    A_k\ge t_k+p-1-o(1),\qquad
    D_k=1-t_k+o(1),\qquad
    G_k=t_k-p-m\varepsilon+o(1).
\]
Since $\liminf\limits_{k\to\infty} t_k\ge\beta=1-\varepsilon$,
$\liminf\limits_{k\to\infty} A_k\ge p-\varepsilon>0$. We have
\begin{align*}
    \liminf_{k\to\infty}\left(G_k-\frac{D_k^2}{A_k}\right)
    \ge 1-m\varepsilon-p-\frac{p\varepsilon}{p-\varepsilon}>\varepsilon.
\end{align*}
The final strict inequality follows from the choice of $\varepsilon$.
For every $\gamma  \in (0, \varepsilon)$ and all sufficiently large $k$,
\begin{equation}\label{eq:thm42-global-descent}
    V_{k+1}^{\mathrm{glob}}-V_k^{\mathrm{glob}}
    \le-\frac{\mu s_k\alpha_k\theta_k(1-\beta)}{\theta_k-\beta}\|z_k-x^\star\|^2
       -\gamma s_k^2W_k\|g_{k-1}\|^2.
\end{equation}

Fix $0<\eta<m\underline\alpha(1-\beta)=m\underline\alpha \varepsilon$ and choose $\eta/m<\gamma<\varepsilon$.
For all sufficiently large $k$, we have
\begin{equation*}
\begin{aligned}
    \gamma s_k^2W_k\|g_{k-1}\|^2
    &\ge2\eta s_kW_k\Delta_{k-1}\\
    &\ge\eta\left[-\frac{s_kW_k}{2L}\|g_{k-1}\|^2
       +2s_kW_k(1-\mu s_{k-1}\varepsilon)\Delta_{k-1}\right].
\end{aligned}
\end{equation*}
Together with \eqref{eq:xi_k_lower}, this bounds the right-hand side of \eqref{eq:thm42-global-descent} by $-\eta V_k^{\mathrm{glob}}$ and proves \eqref{eq:thm42-global-contraction}.
To obtain the lower bound, substitute $\|g_{k-1}\|^2\le2L\Delta_{k-1}$ into the negative gradient term in $V_k^{\mathrm{glob}}$, giving
\[
    V_k^{\mathrm{glob}}
    \ge\|z_k-x^\star\|^2
       +s_kW_k(1-2\mu s_{k-1}\varepsilon)\Delta_{k-1}.
\]
Since $\mu s_{k-1}\to m\le1$ and $\varepsilon<1/8$, the coefficient in parentheses is at least $3/4$ for all sufficiently large $k$. 
This proves \eqref{eq:thm42-global-lower}.
\end{proof}

At least one of the two parts applies: if $\bar q<1$, part~(i) applies; if $\bar q=1$, then $sL\ge\bar q=1>1/2$, so part~(ii) applies. 
When both apply, either energy may be used, and a larger admissible $\eta$ gives a smaller geometric multiplier $1-\eta$.
For the chosen energy $V_k$, increase $N$ so that both its contraction and lower bound hold for every $k\ge N$. 
Iterating the contraction gives $V_{k+1}\le V_N(1-\eta)^{k+1-N}$. 
The energy lower bound further gives $V_{k+1}\ge c_0s_{k+1}W_{k+1}\Delta_k$ for some $c_0>0$.
Since $s_{k+1}$ is bounded away from zero, absorbing fixed factors into $C$ yields
\begin{equation}\label{eq:general-error-rate}
    \Delta_k\le\frac{C}{W_{k+1}}(1-\eta)^{k-N},
    \qquad k\ge N.
\end{equation}
By \textup{(S1)}, $W_k\to\infty$, so this estimate implies $R$-linear convergence of $\Delta_k$. 
Strong convexity gives $\norm{x_k-x^\star}^2\le2\Delta_k/\mu$, establishing $R$-linear convergence of the squared distance as well.

\subsection{A smooth fixed-step FISTA/NAG-type family}

Let $0<s\le1/L$, $r>0$, and $x_0=y_0$. Consider
\begin{equation}\label{eq:smooth-fista-nag}
\begin{cases}
    y_{k+1}=x_k-s\nabla f(x_k),\\[1mm]
    x_{k+1}=y_{k+1}+\dfrac{k}{k+r+1}(y_{k+1}-y_k).
\end{cases}
\end{equation}
This smooth FISTA/NAG-type iteration uses a prescribed schedule for the acceleration coefficients. 
The parameter $r$ is fixed and does not depend on $\mu$.

\begin{corollary}[Fixed-step FISTA/NAG-type convergence]
\label{cor:smooth-fista-general}
Let $f\in\mathcal S^1_{\mu,L}$, and let $\{x_k,y_k\}$ be generated by \eqref{eq:smooth-fista-nag}. Set $\kappa_f:=L/\mu$. For every $\eta>0$ satisfying at least one of
\begin{align}
    \eta&<\mu s(1-sL)=\frac{sL(1-sL)}{\kappa_f},
       &&sL<1,\label{eq:fista-local-rate-condition}\\
    \eta&<\frac{1-\dfrac1{2sL}}{4\kappa_f+1},
       &&sL>\frac12,\label{eq:fista-global-rate-condition}
\end{align}
there are $N\ge1$ and $C>0$ such that
\begin{equation}\label{eq:fista-rate}
    f(x_k)-f_\star\le\frac{C}{k^2}(1-\eta)^{k-N}, \qquad k\ge N.
\end{equation}
The squared distance satisfies the same order:
\[
    \|x_k-x^\star\|^2 =O\!\left(k^{-2}(1-\eta)^k\right).
\]
For example, admissible choices are
\begin{equation}\label{44}
    \eta_{s}:=\begin{cases}
        \frac1{4\kappa_f+1} & \quad \mbox{when } s=1/(2L),\\
        \frac1{8\kappa_f+3} & \quad \mbox{when } s=1/L.
    \end{cases}
\end{equation}
\end{corollary}

\begin{proof}
Set $s_k\equiv s$, $\alpha_k\equiv1$, and $\theta_k:=1+k/r$.
With $z_k:=\theta_kx_k-(\theta_k-1)y_k$, we have $z_0=x_0$ and
\[
    z_{k+1}=z_k-s\theta_kg_k, \qquad
    x_{k+1}=\left(1-\frac1{\theta_{k+1}}\right)y_{k+1}
       +\frac1{\theta_{k+1}}z_{k+1}.
\]
Thus, the iterates satisfy \eqref{eq:scheme}. Moreover,
\[
    W_k=\frac{k(k+r)}{r^2}=\Theta(k^2), \qquad\frac{W_{k+1}}{W_k}\to1, \qquad\frac{H_k}{W_k}=\frac{k+r}{k}\to1.
\]
Since $L_k\le L$, Assumption~\ref{ass:general-structure} holds with $m=\mu s$, $\underline\alpha=1$, and $\bar q\le sL\le1$.

If \eqref{eq:fista-local-rate-condition} holds, choose \(\beta\in\left(sL,1-{\eta}/({\mu s})\right)\).
Then $\bar q\le sL<\beta$ and $\eta<\mu s(1-\beta)$, so Theorem~\ref{thm:general-linear}\textup{(i)} applies.
If \eqref{eq:fista-global-rate-condition} holds, set \(p=1/(2Ls)\) and \(m=\mu s\). 
The parameters in part~\textup{(ii)} give
\begin{equation*}
    m\varepsilon
       =\frac{mp(1-p)}{2+pm}
       =\frac{1-p}{4\kappa_f+1},
\end{equation*}
because $pm=1/(2\kappa_f)$. 
Thus part~\textup{(ii)} applies as well.
In either case, \eqref{eq:general-error-rate} and  $W_{k+1}=\Theta(k^2)$ give \eqref{eq:fista-rate} and the distance estimate.

At $s=1/(2L)$, the local condition admits every $\eta<1/(4\kappa_f)$. 
At $s=1/L$, the global condition admits every $\eta<1/(8\kappa_f+2)$. 
The choices in \eqref{44} lie strictly inside these intervals.
\end{proof}

\begin{remark}[Comparison with the Li--Shi--Yuan bound]
For $0<s<1/L$, let $a:=\mu s(1-sL)>0$.
Li, Shi, and Yuan~\cite{li2024linear} obtain a bound with exponential factor $(1+a/4)^{-k}$. 
The condition \eqref{eq:fista-local-rate-condition} permits $\eta=a/(1+a)$, for which $1-\eta=(1+a)^{-1}$. The corollary therefore gives
\[
    f(x_k)-f_\star=O\!\left(k^{-2}(1+a)^{-k}\right).
\]
This eventual bound improves the exponential factor in that bound while retaining its $k^{-2}$ prefactor. 
The endpoint $s=1/L$ is covered separately by~\eqref{eq:fista-global-rate-condition}.
\end{remark}

\section{Linear convergence without stepsize convergence} \label{sec:nonconvergent}
We now consider stepsizes that need not converge. 

\subsection{Lyapunov contraction and error bounds}
\label{subsec:NC-product}

Suppose that, for some integer $k_0\ge1$, the update and curvature identities required by Lemma~\ref{lem:onestep} hold for every $k\ge k_0$. 
At these indices, assume $s_k>0$, $\alpha_k>0$, $\theta_k>1$, and $\mu\le L_k\le L$.
Define the realized stepsize ratio \(u_k:={s_k}/{s_{k-1}}>0\).
Choose analysis parameters $\beta_k\in(0,1)$ so that
\begin{equation}\label{eq:NC-beta}
    1-\mu s_{k-1}(1-\beta_k)>0.
\end{equation}
Use $F_k$ and $R_k$ from \eqref{eq:common-quadratic-F} and~\eqref{eq:common-quadratic-residual}.
For a nonnegative weight sequence $\omega_k$, take
\begin{equation}\label{eq:NC-energy-choice}
    B_k=\omega_ks_{k-1}^2W_k,\qquad C_k=F_k
\end{equation}
in \eqref{eq:general-energy}, so $\mathcal E_k\ge0$. 
Substitution into
\eqref{eq:common-quadratic-coefficients} gives
\begin{equation}\label{eq:NC-ADG}
\begin{aligned}
    A_k&:=\frac{\beta_k}{s_kL_k}
       -\omega_{k+1}\frac{W_{k+1}}{W_k}
       -\frac{\alpha_k\theta_k}{\theta_k-1}-R_k,\\
    D_k&:=u_k^{-1}-\frac{\beta_k}{s_kL_k},\\
    G_k&:=\frac{\beta_k}{s_kL_k}
       +[\omega_k-\mu s_k(1-\beta_k)]u_k^{-2}.
\end{aligned}
\end{equation}
Here $s_{k-1}=s_ku_k^{-1}$ and $H_k/W_k=\alpha_k\theta_k/(\theta_k-1)$. 
Define the determinant of the coefficient matrix of $Q_k$ by
\begin{equation*}
    \mathfrak m_k
    :=\det\begin{pmatrix}A_k&D_k\\D_k&G_k\end{pmatrix}
    =A_kG_k-D_k^2.
\end{equation*}

\begin{assumption}[Positivity conditions]
\label{ass:NC-margin}
There exists an index $K\ge k_0$ such that, for every $k\ge K$,
\begin{equation*}
    G_k>0,\qquad \mathfrak m_k>0.
\end{equation*}
\end{assumption}

\begin{theorem}[Contraction and error bounds]
\label{thm:NC-determinant}
Let $f\in\mathcal S^1_{\mu,L}$, and suppose that the update and curvature identities required by Lemma~\ref{lem:onestep} hold for every $k\ge K$, where $K\ge k_0$.
For these indices, let $\beta_k\in(0,1)$ and $\omega_k\ge0$, and suppose that \eqref{eq:NC-beta} is satisfied.
Assume also that Assumption~\ref{ass:NC-margin} holds with the same index $K$. Use the energy \eqref{eq:general-energy} with the weights \eqref{eq:NC-energy-choice}.
Define
\begin{equation}\label{eq:NC-rho}
    \rho_k:=\min\left\{
       \frac{\mu s_k\alpha_k\theta_k(1-\beta_k)}{\theta_k-\beta_k},
       \frac{\mathfrak m_k}
         {A_k[\omega_ku_k^{-2}+(\mu s_k)^{-1}]}
    \right\}.
\end{equation}
Then $0<\rho_k<1/4$ and
\begin{equation*}
    \mathcal E_{k+1}\le(1-\rho_k)\mathcal E_k,\qquad k\ge K.
\end{equation*}
Consequently, for $k\ge K$,
\begin{equation}\label{eq:NC-product-rate}
\begin{aligned}
    \Delta_k&\le\frac{\mathcal E_K}{F_{k+1}}
                        \prod_{j=K}^{k}(1-\rho_j),\qquad
    \norm{x_k-x^\star}^2\le\frac{2\mathcal E_K}{\mu F_{k+1}}
                        \prod_{j=K}^{k}(1-\rho_j).
\end{aligned}
\end{equation}
Both bounds remain valid when the product is replaced by $\exp(-\sum_{j=K}^{k}\rho_j)$.
If $\sum_{j=K}^{\infty}\rho_j=\infty$, then $\mathcal E_k\to0$ and $z_k\to x^\star$. 
If also $F_\star:=\inf_{k\ge K}F_k>0$, then $\Delta_k\to0$ and $x_k\to x^\star$. 
If {$F_\star>0$} and $\underline\rho:=\inf_{k\ge K}\rho_k>0$, the objective error and squared distance converge $R$-linearly, with geometric multiplier at most $1-\underline\rho$.
\end{theorem}

\begin{proof}
Since $G_k>0$ and $\mathfrak m_k>0$, we have $A_k>0$.
Completing the square gives
\begin{equation*}
    Q_k\ge\frac{\mathfrak m_k}{A_k}\norm{g_{k-1}}^2.
\end{equation*}
Together with \eqref{eq:common-quadratic-z-descent}, this yields
\begin{equation}\label{eq:NC-basic-descent}
\begin{aligned}
    \mathcal E_{k+1}-\mathcal E_k\le{}&
       -\frac{\mu s_k\alpha_k\theta_k(1-\beta_k)}{\theta_k-\beta_k}
             \norm{z_k-x^\star}^2-\frac{\mathfrak m_k}{A_k}s_k^2W_k\norm{g_{k-1}}^2.
\end{aligned}
\end{equation}
Positivity of $A_k$, nonnegativity of $\omega_{k+1}$ and $R_k$, and
\eqref{eq:NC-ADG} imply
\[
    \mu s_k\frac{\alpha_k\theta_k}{\theta_k-1}
       \le s_kL_k\frac{\alpha_k\theta_k}{\theta_k-1}<\beta_k.
\]
Hence the first entry in the minimum defining $\rho_k$ is strictly smaller than $\beta_k(1-\beta_k)\le1/4$; both entries are positive.

The gradient bound $\Delta_{k-1}\le\norm{g_{k-1}}^2/(2\mu)$ and the positivity in \eqref{eq:NC-beta} give
\begin{equation*}
\begin{aligned}
    B_k\norm{g_{k-1}}^2+F_k\Delta_{k-1}
       &\le s_k^2W_k\left[
          \omega_ku_k^{-2}
          +\frac{1-\mu s_{k-1}(1-\beta_k)}{\mu s_k}
          \right]\norm{g_{k-1}}^2\\
       &\le s_k^2W_k[\omega_ku_k^{-2}+(\mu s_k)^{-1}]
                         \norm{g_{k-1}}^2.
\end{aligned}
\end{equation*}
The definition of $\rho_k$ therefore bounds the right-hand side of~\eqref{eq:NC-basic-descent} by $-\rho_k\mathcal E_k$.
The iteration proves the product bound for $\mathcal E_{k+1}$.
Since $\mathcal E_{k+1}\ge F_{k+1}\Delta_k$ and $\mu\norm{x_k-x^\star}^2/2\le\Delta_k$, we obtain~\eqref{eq:NC-product-rate}; $1-t\le e^{-t}$ gives the exponential replacement. 
The convergence statements follow from these bounds and $\mathcal E_k\ge\norm{z_k-x^\star}^2$.
\end{proof}

\subsection{Explicit conditions for Lyapunov contraction}
\label{subsec:NC-sufficient}

We give sufficient conditions on local curvature, stepsize ratios, and coefficient sequences that allow the analysis parameters $\beta_k$ and $\omega_k$ to be chosen so that~\eqref{eq:NC-beta} and Assumption~\ref{ass:NC-margin} hold.
For $0<\alpha_k<2$, define
\begin{equation*}
    b_k:=\min\{1,2-\alpha_k\},\qquad
    \phi_k:=\frac{b_k}{1+\alpha_kb_k+b_k^2}.
\end{equation*}
Then $0<b_k\le1$, $\alpha_k+b_k\le2$, and $0<\phi_k<1/2$.

\begin{assumption}[Curvature and coefficient conditions]
\label{ass:NC-sufficient}
There exist an index $K\ge k_0$ and constants $\underline s>0$ and $\underline\alpha>0$ such that, for every $k\ge K$, the update and curvature identities required by Lemma~\ref{lem:onestep} hold and the following conditions are satisfied.
\begin{enumerate}
    \item The stepsize and coefficients satisfy
    \begin{equation*}
        s_k\ge\underline s,\qquad
        \underline\alpha\le\alpha_k<2,\qquad\theta_k>1.
    \end{equation*}
    \item The realized stepsize ratio satisfies
    \begin{equation*}
        0<u_k\le U_k,\qquad
        U_k:=\frac{\alpha_{k-1}\theta_{k-1}^2}{W_k}.
    \end{equation*}
    \item Define 
    \begin{equation*}
        v_{0,k}:=\left|\frac{W_{k+1}}{W_k}-1\right|
                +|\alpha_{k+1}-\alpha_k|
                +\frac{\alpha_k}{\theta_k-1}.
    \end{equation*}
    Require
    \begin{equation}\label{eq:NC-P4}
    \begin{aligned}
        \chi_k:= {}&b_k-[1+b_k(\alpha_k+b_k+v_{0,k})]s_kL_k-v_{0,k}U_k^2-2U_k[U_k-1]_+>0.
    \end{aligned}
    \end{equation}
\end{enumerate}
\end{assumption}

Condition \textup{(E2)} limits stepsize growth.
Condition \textup{(E3)} implies that, for every $k\ge K$,
\begin{equation}\label{eq:NC-P2}
    0<s_kL_k<\phi_k<\frac12,\qquad U_k<\frac32.
\end{equation}
Indeed, $\chi_k>0$ gives $s_kL_k<\phi_k$ and $2U_k[U_k-1]_+<b_k\le1$, which excludes $U_k\ge3/2$.

\begin{lemma}[Explicit choice of analysis parameters]
\label{lem:NC-sufficient}
Suppose Assumption~\ref{ass:NC-sufficient} holds with index $K$.
For every $k\ge K$, define the auxiliary parameter
\begin{equation}\label{eq:NC-tau-certificate}
    \tau_k:=\frac{2(U_k+b_k)(U_k+s_kL_k)}{\chi_k}
\end{equation}
and set
\begin{equation}\label{eq:NC-certificate-choice}
    \omega_k:=b_k,\qquad \beta_k:=\frac{\tau_k}{\tau_k+u_k}.
\end{equation}
Then $0<\beta_k<1$ and $1-\mu s_{k-1}(1-\beta_k)>1/2$ for every $k\ge K$. 
Moreover, $A_k>0$, and
\begin{equation}\label{eq:NC-Schur-lower}
    \mathfrak m_k\ge\frac{\beta_k\chi_k}{2s_kL_ku_k^2}>0,
    \qquad
    \frac{\mathfrak m_k}{A_k}\ge\frac{\chi_k}{2U_k^2}.
\end{equation}
Thus \eqref{eq:NC-beta} and Assumption~\ref{ass:NC-margin} hold with the same index $K$.
\end{lemma} 

\begin{proof}
\emph{Bounds on $\beta_k$ and $R_k$.}
The definition of $\chi_k$ gives
\begin{equation}\label{eq:NC-proof-slack-bound}
    0<\chi_k<b_k,\qquad
    [1+b_k(\alpha_k+b_k+v_{0,k})]s_kL_k<b_k.
\end{equation}
In particular, $s_kL_k<b_k$. 
Since $U_k>0$ and $\chi_k<b_k$,
\[
    \tau_k
    >\frac{2b_ks_kL_k}{\chi_k}
    >2s_kL_k\ge2\mu s_k.
\]
In particular, $\tau_k>0$, so $u_k>0$ gives $0<\beta_k<1$.
Hence
\[
    \mu s_{k-1}(1-\beta_k)
       =\frac{\mu s_k}{\tau_k+u_k}<\frac12,
\]
which proves the asserted bound on $1-\mu s_{k-1}(1-\beta_k)$.
Using \textup{(E2)} at index $k+1$, we have
\begin{align*}
    F_{k+1}&\le2s_{k+1}W_{k+1}
             \le2s_k\alpha_k\theta_k^2\\
       &=F_k+2s_k\alpha_k\theta_k
                   +\frac{2\mu s_k^2W_k}{\tau_k+u_k}.
\end{align*}
Taking the positive part and dividing by $2\mu s_k^2W_k>0$ gives
\begin{equation}\label{eq:NC-proof-residual-bound}
    0\le R_k\le\frac1{\tau_k+u_k}.
\end{equation}

\emph{Bounds on the coefficients.}
Fix $k\ge K$ and abbreviate $t=s_kL_k$, $\alpha=\alpha_k$, $b=b_k$, $u=u_k$, $U=U_k$, $v_0=v_{0,k}$, $\chi=\chi_k$, and $\tau=\tau_k$.
The map $\alpha\mapsto\min\{1,2-\alpha\}$ is one-Lipschitz.
Writing $r=W_{k+1}/W_k$, we have $0<b_{k+1}\le1$ and
\[
    b_{k+1}r
    \le b_{k+1}+|r-1|
    \le b+|\alpha_{k+1}-\alpha_k|+|r-1|.
\]
Also, $\alpha_k\theta_k/(\theta_k-1)=\alpha+\alpha/(\theta_k-1)$.
Thus
\[
    b_{k+1}\frac{W_{k+1}}{W_k}
       +\frac{\alpha_k\theta_k}{\theta_k-1}
       \le\alpha+b+v_0.
\]
Substituting $\omega_k=b_k$ into \eqref{eq:NC-ADG} and using~\eqref{eq:NC-proof-residual-bound}, we obtain
\begin{equation}\label{eq:NC-proof-coefficient-bounds}
\begin{gathered}
    0<\frac{\beta_k}{t}-A_k
       \le\alpha+b+v_0+\frac1{\tau+u},\\
    0<b-\frac{tu}{\tau+u}
       \le b-\mu s_k(1-\beta_k)\le b,\\
    0<G_k\le\frac{\beta_k}{t}+\frac b{u^2}.
\end{gathered}
\end{equation}
The first expression is positive because $\beta_k/t-A_k=b_{k+1}W_{k+1}/W_k+\alpha_k\theta_k/(\theta_k-1)+R_k>0$. The second line uses $\mu s_k\le t<b$ and $u/(\tau+u)<1$. Note it also proves $G_k>0$.

\emph{Lower bound on the determinant.}
Expanding $\mathfrak m_k=A_kG_k-D_k^2$ gives
\[
    \mathfrak m_k
       =\frac{\beta_k}{tu^2}
           [b-\mu s_k(1-\beta_k)+2u]
         -\frac1{u^2}
         -\left(\frac{\beta_k}{t}-A_k\right)G_k.
\]
Use the lower bound for $b-\mu s_k(1-\beta_k)$ and the upper bounds for the two positive factors in the last product from~\eqref{eq:NC-proof-coefficient-bounds}.
Substituting $1/\beta_k=1+u/\tau$ and bounding $1/(\tau+u)\le1/\tau$ give
\begin{equation}\label{eq:NC-proof-margin-lower}
\begin{aligned}
    \mathfrak m_k\ge\frac{\beta_k}{tu^2}\Biggl[
       &b-[1+b(\alpha+b+v_0)]t+2u-(\alpha+b+v_0)u^2\\
       &-\frac{bt+[2+b(\alpha+b+v_0)]tu+u^2}{\tau}
    \Biggr].
\end{aligned}
\end{equation}
Because $\alpha+b\le2$ and $0<u\le U$,
\[
    2u-(\alpha+b)u^2\ge-2U[U-1]_+,
    \qquad -v_0u^2\ge-v_0U^2.
\]
The terms in brackets that do not contain $\tau$ are therefore at least $\chi$. The bound~\eqref{eq:NC-proof-slack-bound} also gives $[2+b(\alpha+b+v_0)]t<b+t$. 
Consequently,
\[
    bt+[2+b(\alpha+b+v_0)]tu+u^2
       \le bt+(b+t)U+U^2=(U+b)(U+t).
\]
By the choice of $\tau$, the last fraction in~\eqref{eq:NC-proof-margin-lower} is at most $\chi/2$. 
This proves the first bound in \eqref{eq:NC-Schur-lower}.

Since $G_k>0$ and $\mathfrak m_k>0$, we have $A_k>0$.
The first line of \eqref{eq:NC-proof-coefficient-bounds} also gives $A_k<\beta_k/t$. 
Hence
\[
    \frac{\mathfrak m_k}{A_k}
    \ge\frac{t}{\beta_k}\mathfrak m_k
    \ge\frac{\chi}{2u^2}
    \ge\frac{\chi}{2U^2}.
\]
This proves \eqref{eq:NC-Schur-lower}. Since $k\ge K$ was arbitrary, all conclusions hold for every $k\ge K$.
\end{proof}

\begin{corollary}[Explicit error bounds]
\label{cor:NC-sufficient}
Suppose Assumption~\ref{ass:NC-sufficient} holds with index $K$.
Define $\tau_k$ by \eqref{eq:NC-tau-certificate} and choose $\beta_k,\omega_k$ according to \eqref{eq:NC-certificate-choice}.
For every $k\ge K$, define
\begin{equation}\label{eq:NC-certificate-rate}
    \widehat\rho_k:=\min\left\{
       \frac{\mu s_k\alpha_k\theta_ku_k}
            {(\tau_k+u_k)(\theta_k-1)+u_k},
       \frac{\chi_k}
            {2U_k^2[b_ku_k^{-2}+(\mu s_k)^{-1}]}
    \right\}.
\end{equation}
Then $0<\widehat\rho_k\le\rho_k<1/4$, where $\rho_k$ is defined
in \eqref{eq:NC-rho}, and
\[
    \mathcal E_{k+1}
    \le(1-\widehat\rho_k)\mathcal E_k,
    \qquad k\ge K.
\]
Consequently,
\begin{equation}\label{eq:NC-cor-order}
    \Delta_k,\ \norm{x_k-x^\star}^2
    =O\!\left(
       \frac1{W_{k+1}}
       \exp\!\left[-\sum_{j=K}^{k}\widehat\rho_j\right]
    \right).
\end{equation}
Moreover, $\inf_{k\ge K}F_k>0$.
If $\sum_{j=K}^{\infty}\widehat\rho_j=\infty$, then $\mathcal E_k\to0$, $z_k\to x^\star$, $\Delta_k\to0$, and $x_k\to x^\star$.
If \(\widehat\rho_*:=\inf_{k\ge K}\widehat\rho_k>0\), then the objective error and squared distance converge $R$-linearly, with geometric multiplier at most $1-\widehat\rho_*$.
\end{corollary}

\begin{proof}
The first term in \eqref{eq:NC-certificate-rate} is the first term in~\eqref{eq:NC-rho}, after substituting $1-\beta_k=u_k/(\tau_k+u_k)$. 
With $\omega_k=b_k$, the Schur bound~\eqref{eq:NC-Schur-lower} supplies the second term. 
Hence $0<\widehat\rho_k\le\rho_k<1/4$, and Theorem~\ref{thm:NC-determinant} implies $\mathcal E_{k+1}\le(1-\widehat\rho_k)\mathcal E_k$.

By Lemma~\ref{lem:NC-sufficient}, $\mu s_k/(\tau_k+u_k)<1/2$. Therefore,
\[
    F_k
    =2s_kW_k\left(1-\frac{\mu s_k}{\tau_k+u_k}\right)
    >s_kW_k
    \ge\underline s\,W_k.
\]
This proves \eqref{eq:NC-cor-order}. 
For $k\ge K+1$, \textup{(E1)} gives $\alpha_{k-1}\ge\underline\alpha$ and $\theta_{k-1}>1$, while \eqref{eq:NC-P2} gives $U_k<3/2$. 
Hence
\(W_k
    ={\alpha_{k-1}\theta_{k-1}^2}/{U_k}
    >2\underline\alpha/3\). 
Since $W_K>0$, the lower bound on $F_k$ yields
\(\inf_{k\ge K}F_k
    \ge\underline s\min\left\{W_K,2\underline\alpha/3\right\}
    >0\).
The convergence statements now follow from Theorem~\ref{thm:NC-determinant}, since $\rho_k\ge\widehat\rho_k$ and $F_k$ is bounded away from zero.
\end{proof}

\subsection{An explicit starting index and linear convergence}

The following theorem gives an explicit index $K$ after which $\chi_k$ is bounded away from zero. 
Direct verification of $\chi_k>0$ may yield a smaller starting index.

\begin{theorem}[Explicit starting index for the sufficient conditions]
\label{thm:unified-tail-verification}
Let $f\in\mathcal S^1_{\mu,L}$. 
Suppose that the update and curvature identities required by Lemma~\ref{lem:onestep}, together with \textup{(E1)}--\textup{(E2)}, hold for every $k\ge k_0$, where $k_0\ge1$.
Assume there exist constants $r_\theta>0$, $C_0\ge0$, $C_U\ge0$, $\overline U>0$, and $\delta_{\rm alg}>0$ such that, for every $k\ge k_0$,
\begin{equation}\label{eq:unified-direct-tail-bounds}
    v_{0,k}\le\frac{C_0}{k^{r_\theta}},\qquad
    0<U_k\le\overline U,\qquad
    [U_k-1]_+\le\frac{C_U}{k^{r_\theta}},
\end{equation}
and
\begin{equation}\label{eq:unified-delta-margin}
    b_k-(1+\alpha_kb_k+b_k^2)s_kL_k
    \ge\delta_{\rm alg}.
\end{equation}
Set
\begin{equation}\label{eq:unified-entrance-index}
    K:=\max\left\{k_0,\ 1+\left\lfloor
       \left[\frac{C_0(\overline U^2+1/2)+2\overline U C_U}
                       {\delta_{\rm alg}}\right]^{1/r_\theta}
       \right\rfloor\right\}.
\end{equation}
Then Assumption~\ref{ass:NC-sufficient} holds with index $K$.
Moreover,
\begin{align*}
    \inf_{k\ge K}\chi_k
    &>0,\qquad \liminf_{k\to\infty}\chi_k
    \ge\delta_{\rm alg}.
\end{align*}
\end{theorem}

\begin{proof}
Since $\delta_{\rm alg}>0$,~\eqref{eq:unified-delta-margin} implies $s_kL_k<\phi_k<1/2$.
Set $C_\chi:=C_0(\overline U^2+1/2)+2\overline U C_U$.
Using $b_k\le1$ and \eqref{eq:unified-direct-tail-bounds}, we obtain
\begin{align*}
    \chi_k&\ge\delta_{\rm alg}
          -v_{0,k}(U_k^2+b_ks_kL_k)-2U_k[U_k-1]_+\\
       &\ge\delta_{\rm alg}
          -\frac{C_\chi}{k^{r_\theta}}
        \ge\delta_{\rm alg}-\frac{C_\chi}{K^{r_\theta}}>0,
        \qquad k\ge K.
\end{align*}
The last inequality follows from $K>(C_\chi/\delta_{\rm alg})^{1/r_\theta}$. Thus
$\inf_{k\ge K}\chi_k\ge\delta_{\rm alg}-C_\chi/K^{r_\theta}>0$.
Thus \textup{(E3)} holds for every $k\ge K$.
Since \textup{(E1)}--\textup{(E2)} already hold at these indices, Assumption~\ref{ass:NC-sufficient} holds with index $K$.
Finally, $1/k^{r_\theta}\to0$ gives $\liminf_{k\to\infty}\chi_k\ge\delta_{\rm alg}$.
\end{proof}

\begin{corollary}[Eventual linear convergence]
\label{cor:unified-uniform-rate}
Suppose the hypotheses of Theorem~\ref{thm:unified-tail-verification} hold, and let $K$ be given by \eqref{eq:unified-entrance-index}.
For $k\ge K$, define $\tau_k$ by~\eqref{eq:NC-tau-certificate} and choose $\beta_k,\omega_k$ according to \eqref{eq:NC-certificate-choice}.
Set $\sigma_G:=L\underline s$ and $\kappa_f:=L/\mu$, and define
\begin{equation*}
\begin{aligned}
    c_k:=\min\Biggl\{&
       \frac{\sigma_G\alpha_k\theta_ku_k}
            {(\tau_k+u_k)(\theta_k-1)+u_k},\quad
       \frac{\sigma_G\chi_k}
            {2U_k^2(1+b_ks_kL_ku_k^{-2})}
    \Biggr\}.
\end{aligned}
\end{equation*}
Then, with $\widehat\rho_k$ defined in
\eqref{eq:NC-certificate-rate},
\begin{equation}\label{eq:unified-rho-lower-bound}
    \widehat\rho_k\ge\frac{c_k}{\kappa_f}>0,
    \qquad k\ge K.
\end{equation}
Consequently,
\begin{equation}\label{eq:unified-product-rate}
    \Delta_k,\ \norm{x_k-x^\star}^2
    =O\!\left(
       \frac1{W_{k+1}}
       \exp\!\left[-\frac1{\kappa_f}\sum_{j=K}^{k}c_j\right]
    \right).
\end{equation}
Moreover, \(c:=\inf_{k\ge K}c_k>0\) and \(0< c/{\kappa_f}<\frac14\).
Thus, the objective error and the squared distance converge $R$-linearly, with a geometric multiplier at most $1-c/\kappa_f$.
If, in addition, $W_k=\Theta(k^{2r_\theta})$, then
\begin{equation}\label{eq:unified-linear-rate}
    \Delta_k,\ \norm{x_k-x^\star}^2
    =O\!\left(
       k^{-2r_\theta}
       \left(1-\frac c{\kappa_f}\right)^{k-K}
    \right).
\end{equation}
\end{corollary}

\begin{proof}
For the first term in \eqref{eq:NC-certificate-rate}, use $\mu s_k\ge\mu\underline s=\sigma_G/\kappa_f$.
For the second term, $\mu\le L_k$ and $s_k\ge\underline s$ give
\[
\begin{aligned}
    \frac{\chi_k}
         {2U_k^2[b_ku_k^{-2}+(\mu s_k)^{-1}]}
    &=\frac{\mu s_k\chi_k}
            {2U_k^2(1+b_k\mu s_ku_k^{-2})}\ge \frac1{\kappa_f}
      \frac{\sigma_G\chi_k}
           {2U_k^2(1+b_ks_kL_ku_k^{-2})}.
\end{aligned}
\]
Comparing the two terms proves~\eqref{eq:unified-rho-lower-bound}. Together with Corollary~\ref{cor:NC-sufficient}, it yields \eqref{eq:unified-product-rate}.

We next establish a positive lower bound on $c_k$ independent of $k$.
By Theorem~\ref{thm:unified-tail-verification},
$\chi_*:=\inf_{k\ge K}\chi_k>0$.
Condition \textup{(E1)} and \eqref{eq:NC-P2} give
\[
    \underline s\le s_k<\frac1{2\mu},
    \qquad k\ge K.
\]
Thus $u_k=s_k/s_{k-1}>2\mu\underline s$ for $k\ge K+1$.
Including the first value $u_K>0$, we obtain
\(u_*:=\inf_{k\ge K}u_k
    \ge\min\{u_K,2\mu\underline s\}>0\).
Also, $u_k\le U_k<3/2$, $b_k\le1$, and $s_kL_k<1/2$.
Consequently,
\[
    \tau_k
    =\frac{2(U_k+b_k)(U_k+s_kL_k)}{\chi_k}
    \le\frac{2(3/2+1)(3/2+1/2)}{\chi_*}
    =\frac{10}{\chi_*}.
\]

Write $T:=10/\chi_*$. Since
\[
    (\tau_k+u_k)(\theta_k-1)+u_k
    =\theta_k(\tau_k+u_k)-\tau_k
    <\theta_k(\tau_k+u_k),
\]
the first term defining $c_k$ is at least
\[
    \frac{\sigma_G\underline\alpha u_k}{\tau_k+u_k}
    \ge\frac{\sigma_G\underline\alpha u_*}{T+3/2}>0.
\]
The second has a positive infimum because its numerator is at least $\sigma_G\chi_*$ and its denominator is at most $2(3/2)^2(1+u_*^{-2}/2)$. 
Hence $c:=\inf_{k\ge K}c_k>0$.
Moreover, $c/\kappa_f\le c_k/\kappa_f\le\widehat\rho_k<1/4$. 
The contraction with $\widehat\rho_k\ge c/\kappa_f$ gives
\[
    \Delta_k
    \le\frac{\mathcal E_K}{F_{k+1}}
           (1-c/\kappa_f)^{k+1-K}
    \le\frac{\mathcal E_K}{\underline s\,W_{k+1}}
           (1-c/\kappa_f)^{k+1-K}.
\]
The positive lower bound on $F_k$ established in Corollary~\ref{cor:NC-sufficient} proves $R$-linear convergence of $\Delta_k$. Strong convexity gives the same conclusion for $\norm{x_k-x^\star}^2$. If $W_k=\Theta(k^{2r_\theta})$, the last display gives \eqref{eq:unified-linear-rate}.
\end{proof}

 \section{Applications and framework-guided design}
\label{sec:framework-consequences}

We first verify the conditions of Section~\ref{sec:nonconvergent} for existing methods, then use them to construct two new families.
For each policy, we prove positive lower bounds on the stepsizes and coupling coefficients, check $u_k\le U_k$, and bound
$b_k-(1+\alpha_kb_k+b_k^2)s_kL_k$ below by a positive constant.
Together with bounds on coefficient changes, these inequalities allow Theorem~\ref{thm:unified-tail-verification} and Corollary~\ref{cor:unified-uniform-rate} to be applied.

\subsection{Verification of existing methods}
\label{subsec:existing-algorithms}

We use the representations in Section~\ref{sec:setting} and the local curvature rule \eqref{eq:Lk-bregman-form}, including its continuation for coincident iterates.
The fixed algorithmic parameters $\beta^{\rm AC}$ and $\beta^{\rm AGM}$ are distinct from the analysis sequence $\beta_k$.

\subsubsection{AdaNAG-G}
\label{subsec:G-consequences}
For AdaNAG-G~\cite{suh2025adanag}, choose coefficient sequences $\alpha_k>0$ and $\theta_k>1$ satisfying
\((\theta_{k+1}-1)^2>\alpha_k\theta_k^2\) for every $k\ge0$.
For $k\ge0$, define
\begin{align*}
    A_k^G&:=\alpha_{k+1}\theta_{k+1}(\theta_{k+1}-1),\\
    B_{k+1}^G&:=\alpha_{k+1}^2\theta_{k+1}^2
       \left(
          \frac{(\theta_{k+1}-1)^2}{\alpha_k\theta_k^2}-1
       \right).
\end{align*}
The coefficient restriction ensures $B_{k+1}^G>0$.
Supply $s_0>0$, $A_{-1}^G=0$, and $B_0^G>0$ for initialization.
For every $k\ge0$, update the stepsize by
\begin{equation}\label{eq:adanagg-stepsize}
\begin{aligned}
    s_{k+1}=\min\Biggl\{&
       \frac{A_{k-1}^G+\alpha_k\theta_k}{A_k^G}s_k,\\[-1mm]
       &\left(
          \frac{A_k^G}{B_k^G}
          +\frac{B_{k+1}^G+\alpha_{k+1}^2\theta_{k+1}^2}{A_k^G}
        \right)^{-1}\frac1{L_{k+1}}
    \Biggr\}.
\end{aligned}
\end{equation}
For $k\ge2$, the second branch gives
\begin{equation}\label{eq:G-general-curvature-certificate}
    s_kL_k\le q_k^G:=
    \left(
       \frac{A_{k-1}^G}{B_{k-1}^G}
       +\frac{B_k^G+\alpha_k^2\theta_k^2}{A_{k-1}^G}
    \right)^{-1}.
\end{equation}
Moreover, $A_{k-1}^G=W_k$ and \(A_{k-2}^G+\alpha_{k-1}\theta_{k-1} =\alpha_{k-1}\theta_{k-1}^2\).
Thus, for $k\ge2$, the stepsize rule can be written as
\[
    s_k=\min\left\{U_ks_{k-1},\frac{q_k^G}{L_k}\right\},
\]
which also gives $u_k\le U_k$.

We first establish linear convergence without requiring a prescribed growth rate for $\theta_k$.

\begin{corollary}[Linear convergence of AdaNAG-G]\label{cor:G-general-coupled}
Let $f\in\mathcal S^1_{\mu,L}$, and consider the AdaNAG-G policy
\eqref{eq:adanagg-stepsize}. Suppose its coefficients satisfy
\begin{equation*}
    0<\alpha_k\le1,\qquad
    \alpha_k\to\bar\alpha\in(0,1),\qquad
    \theta_k\to\infty,\qquad
    \frac{\theta_{k+1}}{\theta_k}\to1,
\end{equation*}
and suppose there is an index $N_0\ge1$ such that
\begin{equation*}
    U_k\ge1,\qquad k\ge N_0.
\end{equation*}
Then there exists a finite index $K_G\ge N_0$ such that Assumption~\ref{ass:NC-sufficient} holds with index $K_G$ and
\begin{equation}\label{eq:G-general-rho-liminf}
    \liminf_{k\to\infty}\widehat\rho_k>0.
\end{equation}
Consequently, $\Delta_k$ and $\norm{x_k-x^\star}^2$ converge $R$-linearly to zero.
\end{corollary}

\begin{proof}
The coefficient limits give
$W_{k+1}/W_k\to1$, $v_{0,k}\to0$, and $U_k\to1$.
They also imply
\[
    q_k^G\longrightarrow q_\infty:=\frac{1-\bar\alpha}{2-\bar\alpha}>0.
\]
Choose $N\ge\max\{N_0,2\}$ large enough that $q_k^G\ge q_\infty/2$ and $\alpha_k\ge\bar\alpha/2$ for every $k\ge N$.
Since $U_k\ge1$ and $L_k\le L$, the stepsize rule gives
\[
    s_k\ge\underline s:=
    \min\left\{s_{N-1},\frac{q_\infty}{2L}\right\}>0,
    \qquad k\ge N.
\]
Together with $\alpha_k\ge\bar\alpha/2$, this verifies~\textup{(E1)}. 
The stepsize rule also gives \textup{(E2)}.

Since $b_k=1$, $s_kL_k\le q_k^G$, $v_{0,k}\to0$, and $U_k\to1$,
\[
    \liminf_{k\to\infty}\chi_k
    \ge1-(2+\bar\alpha)q_\infty
    =\frac{\bar\alpha^2}{2-\bar\alpha}>0.
\]
Thus we can choose $K_G\ge N$ such that
$\inf_{k\ge K_G}\chi_k>0$, establishing \textup{(E3)}.

For $k\ge K_G$, \eqref{eq:NC-P2} gives $s_k<1/(2\mu)$.
Hence
\[
    u_k=\frac{s_k}{s_{k-1}}>2\mu\underline s,
    \qquad k\ge K_G+1.
\]
The formula \eqref{eq:NC-tau-certificate} bounds $\tau_k$ above because $U_k$ is bounded, $s_kL_k<1/2$, and $\inf_{k\ge K_G}\chi_k>0$.
These estimates and the positive lower bounds on $s_k$ and $\alpha_k$ give positive lower limits for both terms in \eqref{eq:NC-certificate-rate}.
This proves \eqref{eq:G-general-rho-liminf}.
Corollary~\ref{cor:NC-sufficient} then gives the stated $R$-linear convergence.
\end{proof}

The two coefficient sequences used in the experiments are, for $k\ge0$,
\begin{equation}\label{eq:G-concrete-policies}
\begin{aligned}
    \text{AdaNAG-G}_{12}:&\quad
       \theta_k=\frac{k+14}{12},\qquad
       \alpha_k=\frac12\left(\frac{k+3}{k+14}\right)^2,\\
    \text{AdaNAG-G}_{1/2}:&\quad
       \theta_k=2\sqrt{k+3},\qquad
       \alpha_k=\frac12.
\end{aligned}
\end{equation}
For both policies,
\[
    U_k>1,\qquad
    \frac1{10}\le q_k^G<\frac13, \qquad k\ge4,
\]
as shown in Appendix~\ref{app:section6-native-g-checks}.
The stepsize rule and $0<\alpha_k\le1/2$ therefore give
\[
    s_k\ge\min\left\{s_3,\frac1{10L}\right\},\qquad 1-(2+\alpha_k)s_kL_k
    \ge1-\frac{2+\alpha_k}{3}
    \ge\frac16,
    \qquad k\ge4.
\]
Thus \eqref{eq:unified-delta-margin} holds with
$\delta_{\rm alg}=1/6$.

Lemma~\ref{lem:unified-coefficient-bounds} supplies the coefficient
bounds required by Theorem~\ref{thm:unified-tail-verification},
with $r_\theta=1$ and $r_\theta=1/2$, respectively.
It also gives $W_k=\Theta(k^{2r_\theta})$.
Corollary~\ref{cor:unified-uniform-rate} therefore yields
\eqref{eq:unified-linear-rate}, with polynomial prefactors
$k^{-2}$ and $k^{-1}$, respectively.
The initialization used in the experiments is specified in
Section~\ref{sec:numerics}.

\subsubsection{Adaptive AC-FGM}
We consider stepsize policy II for smooth unconstrained
AC-FGM~\cite{li2026uniformly}, with
$a\in(0,1]$ and $0<\beta^{\rm AC}\le1-\sqrt6/3$.
Initialize $x_0=z_0=x^{\rm init}$ and choose
$\eta^{\rm init}>0$.
If $g_0=0$, stop at the minimizer. Otherwise, set
\[
    x_1=x^{\rm init}-\eta^{\rm init}g_0,\qquad
    z_1=x^{\rm init},\qquad
    \tau_1^{\rm AC}=0,\qquad \tau_2^{\rm AC}=1.
\]
Define the initial secant estimate and the AC-FGM stepsize $\eta_2$ by
\[
    L_1=\frac{\norm{g_1-g_0}}{\norm{x_1-x_0}},\qquad
    \eta_2=\min\left\{
       (1-\beta^{\rm AC})\eta^{\rm init},\frac1{4L_1}
    \right\}.
\]
Since $x_1\ne x_0$, the secant estimate satisfies $\mu\le L_1\le L$.
The second update, expressed using the unified intermediate point $y_2$, gives
\[
\begin{aligned}
    y_2&=x_1-(1-\beta^{\rm AC})\eta_2g_1,
    &z_2=x^{\rm init}-\beta^{\rm AC}\eta_2g_1,\\
    x_2&=\frac{y_2+z_2}{2},\qquad \theta_2=2.&
\end{aligned}
\]
To express the initialization in the common notation, set
\(
    \theta_1=1,~
    \alpha_1=\beta^{\rm AC}/(1-\beta^{\rm AC}),~s_1=(1-\beta^{\rm AC})\eta_2.
\)
These values give
$y_2=x_1-s_1g_1$ and $z_2=2x_2-y_2$.

For $k\ge2$, compute $L_k$ by
\eqref{eq:Lk-bregman-form} and update
\begin{align*}
    \eta_{k+1}
    &=\min\left\{
       \frac43\eta_k,\,
       \frac{1+\tau_{k-1}^{\rm AC}}{\tau_k^{\rm AC}}\eta_k,\,
       \frac{\tau_k^{\rm AC}}{4L_k}
    \right\},\\
    \tau_{k+1}^{\rm AC}
    &=\tau_k^{\rm AC}+\frac a2
       +\frac{2(1-a)\eta_{k+1}L_k}{\tau_k^{\rm AC}}.
\end{align*}
With the identification
\eqref{eq:acfgm-unified-parameters}, these updates satisfy the
unified iteration.
The initialization supplies the preceding gradient-step and coupling
identities at $k=2$. With $L_2$ computed by the Bregman rule,
the update and curvature identities required by
Lemma~\ref{lem:onestep} therefore hold from that index.

The unified stepsize satisfies
\begin{equation*}
    s_k=\min\left\{
       \frac{4\tau_k^{\rm AC}}{3\tau_{k+1}^{\rm AC}}s_{k-1},\,
       U_ks_{k-1},\,
       \frac{q_k^{\rm AC}}{L_k}
    \right\},\qquad k\ge2,
\end{equation*}
where
\begin{equation*}
    q_k^{\rm AC}:=
       \frac{(1-\beta^{\rm AC})\tau_k^{\rm AC}}
            {4\tau_{k+1}^{\rm AC}},\qquad
    U_k=\frac{1+\tau_{k-1}^{\rm AC}}{\tau_{k+1}^{\rm AC}}.
\end{equation*}
In particular, $u_k\le U_k$ and $s_kL_k\le q_k^{\rm AC}$.
The bound $\eta_{k+1}\le\tau_k^{\rm AC}/(4L_k)$ and the update of $\tau_k^{\rm AC}$ give
\[
    \frac a2\le\tau_{k+1}^{\rm AC}-\tau_k^{\rm AC}\le\frac12,
    \qquad
    \frac{1-\beta^{\rm AC}}6\le q_k^{\rm AC}
       \le\frac{1-\beta^{\rm AC}}4,\qquad k\ge2.
\]

The increment bound gives $U_k\ge1$ for $k\ge3$, since
$\tau_{k+1}^{\rm AC}-\tau_{k-1}^{\rm AC}\le1$.
It also gives
\[
    \tau_k^{\rm AC}\ge1+\frac{a(k-2)}2,\qquad k\ge2.
\]
Thus $N:=2+\lceil1/a\rceil$ ensures $\tau_k^{\rm AC}\ge3/2$ and
$4\tau_k^{\rm AC}/(3\tau_{k+1}^{\rm AC})\ge1$ for every $k\ge N$.
Since $L_k\le L$, induction gives
\[
    s_k\ge\underline s:=\min\left\{
       \min_{1\le j<N}s_j,\frac{1-\beta^{\rm AC}}{6L}\right\}>0,
       \qquad k\ge1.
\]
The initial multiplier $U_2=1/\tau_3^{\rm AC}<1$ is included in the finite prefix.

For $k\ge2$, the coefficients satisfy
\[
    \frac{\beta^{\rm AC}}{2(1-\beta^{\rm AC})}
       \le\alpha_k<\frac{\beta^{\rm AC}}{1-\beta^{\rm AC}}<1,
       \qquad \theta_k>1.
\]
This proves \textup{(E1)}--\textup{(E2)} for every $k\ge2$. Since $b_k=1$,
\begin{equation}\label{eq:acfgm-curvature-margin}
    1-(2+\alpha_k)s_kL_k
       \ge1-(2+\alpha_k)q_k^{\rm AC}
       \ge\frac{2+\beta^{\rm AC}}4>0.
\end{equation}
To estimate coefficient changes, put $h=\beta^{\rm AC}/(1-\beta^{\rm AC})$.
The bounded increments of $\tau_k^{\rm AC}$ give
\[
    \tau_k^{\rm AC}=\Theta(k),\qquad
    \alpha_k=h+O(k^{-1}),\qquad
    W_k=h\tau_k^{\rm AC}\tau_{k+1}^{\rm AC}=\Theta(k^2).
\]
They also give $v_{0,k}=O(k^{-1})$ and $[U_k-1]_+=O(k^{-1})$. 
These bounds and \eqref{eq:acfgm-curvature-margin} verify the hypotheses of Theorem~\ref{thm:unified-tail-verification}, with $r_\theta=1$ and $\delta_{\rm alg}=(2+\beta^{\rm AC})/4$.
Corollary~\ref{cor:unified-uniform-rate} therefore gives~\eqref{eq:unified-linear-rate}.

\subsubsection{AdaAGM}
For AdaAGM~\cite{wang2025adaagm}, choose $0<m<1$, $0<\gamma<2$, $\beta^{\rm AGM}>0$, $\theta_0\ge1$, and $\Omega,~\eta\in[0,1)$ satisfying
\begin{equation}\label{eq:adaagm-initialization}
    \frac{(1+\beta^{\rm AGM})\gamma}{2}
    \le\frac{2(\theta_0-1)+m}{2\theta_0+m}.
\end{equation}
For our verification, exclude $(\Omega,\eta,\beta^{\rm AGM}\gamma)=(0,0,1)$.
This restriction ensures the positive asymptotic margin established below.
In particular, the verified family includes the admissible parameter choices for which Theorem~4.1 of~\cite{wang2025adaagm} gives a positive linear convergence rate.

Initialize $x_0=y_0=z_0=x^{\rm init}$, choose $s_0>0$, and set
\[
    M_{\rm AGM}:=
    \frac{1-\Omega}
         {2+1/[\beta^{\rm AGM}\gamma(1-\eta)]},
    \qquad
    \underline s:=\min\{s_0,M_{\rm AGM}/L\}>0.
\]
The published strongly convex theorem also requires $s_0\ge M_{\rm AGM}/L$.
The local choice $s_0=M_{\rm AGM}/L_0$, with $L_0\in[\mu,L]$, satisfies this requirement without using the global constant $L$.

For $k\ge0$, define
\begin{equation}\label{eq:adaagm-theta}
    \theta_{k+1}
    =\frac{m+\sqrt{m^2+4\theta_k^2}}2,
    \qquad \alpha_k\equiv\gamma,
\end{equation}
and use these coefficients in~\eqref{eq:adaagm-update}.
If $g_0=0$, stop at the minimizer. Otherwise, we have \(x_1\ne x_0\).
The initialization supplies the previous gradient-step and coupling identities required by Lemma~\ref{lem:onestep} from $k=1$.

For $k\ge1$, define
\begin{equation}\label{eq:adaagm-curvature-certificate}
    q_k^{\rm AGM}:=
    \frac{1-\Omega}{
       (1+\beta^{\rm AGM})\gamma
          \dfrac{\theta_k}{\theta_k-1}
       +\dfrac{\theta_k-1}
          {\beta^{\rm AGM}\gamma(1-\eta)(\theta_k-m)}}.
\end{equation}
The native stepsize rule is
\begin{equation}\label{eq:adaagm-stepsize}
\begin{aligned}
    s_{k+1}=\min\Biggl\{&
       U_{k+1}s_k,\,
       \frac{2}{(1+\beta^{\rm AGM})\gamma}
          \left(1-\frac1{\theta_{k+1}}\right)s_k,\,
       \frac{q_{k+1}^{\rm AGM}}{\widehat L_{k+1}}
    \Biggr\},
    ~~k\ge0,
\end{aligned}
\end{equation}
where $\widehat L_k=L_k$ for $x_k\ne x_{k-1}$ and $\widehat L_k=0$ otherwise, with $q_k^{\rm AGM}/0:=+\infty$.

The recurrence gives $\theta_{k+1}>\theta_k+m/2$.
Consequently,~\eqref{eq:adaagm-initialization} implies
\[
    \frac{(1+\beta^{\rm AGM})\gamma}{2}
    \le1-\frac1{\theta_0+m/2}
    <1-\frac1{\theta_1},\quad U_k
    =\frac{\theta_k-m}{\theta_k-1}>1,
\]
Since $\theta_k\ge\theta_1$ for $k\ge1$, both growth multipliers
in~\eqref{eq:adaagm-stepsize} exceed one.
Moreover,
\[
    \frac{(1+\beta^{\rm AGM})\gamma\theta_k}{\theta_k-1}<2,
    \qquad
    0<\frac{\theta_k-1}{\theta_k-m}<1,
    \qquad k\ge1,
\]
so $q_k^{\rm AGM}>M_{\rm AGM}$.
Induction yields \(s_k\ge\underline s\). Together with $\alpha_k=\gamma\in(0,2)$ and $\theta_k>1$, these bounds verify \textup{(E1)}--\textup{(E2)}.

To verify~\textup{(E3)}, write $b_k\equiv b:=\min\{1,2-\gamma\}$ and define
\begin{equation*}
\begin{aligned}
    q_\infty&:=
    \frac{1-\Omega}{
       (1+\beta^{\rm AGM})\gamma
       +1/[\beta^{\rm AGM}\gamma(1-\eta)]},\\
    \delta_{\rm AGM}&:=b-(1+\gamma b+b^2)q_\infty.
\end{aligned}
\end{equation*}
We first establish $\delta_{\rm AGM}>0$.
For $0<\gamma\le1$, the arithmetic--geometric mean inequality gives
\[
    \beta^{\rm AGM}\gamma
    +\frac1{\beta^{\rm AGM}\gamma(1-\eta)}
    \ge\frac2{\sqrt{1-\eta}}\ge2.
\]
Hence $q_\infty\le1/(\gamma+2)$.
Equality requires $(\Omega,\eta,\beta^{\rm AGM}\gamma)=(0,0,1)$, which is excluded.
Therefore, \(q_\infty<1/(\gamma+2)=b/(1+\gamma b+b^2)\).
For $1<\gamma<2$,~\eqref{eq:adaagm-initialization} implies \(0<\beta^{\rm AGM}\gamma<2-\gamma=b<1\).
Since $t+1/t$ is decreasing on $(0,1)$,
\[
    (1+\beta^{\rm AGM})\gamma
    +\frac1{\beta^{\rm AGM}\gamma(1-\eta)}
    >\gamma+b+\frac1b
    =2+\frac1b.
\]
Thus \(q_\infty<{b}/({1+2b})={b}/({1+\gamma b+b^2})\) in this case as well. This proves $\delta_{\rm AGM}>0$.

The recurrence~\eqref{eq:adaagm-theta} gives
\[
    \frac m2<\theta_{k+1}-\theta_k<m,\qquad
    \theta_k=\Theta(k),\qquad W_k=\Theta(k^2),
\]
and \(0<W_{k+1}-W_k
      =\gamma\bigl[m\theta_k
          -(1-m)(\theta_{k+1}-\theta_k)\bigr]
      <m\gamma\theta_k.\)
Consequently, for $k\ge1$,
\[
    v_{0,k}<\frac{m+\gamma}{\theta_k-1}=O(k^{-1}),
    \qquad
    U_k=1+\frac{1-m}{\theta_k-1}=1+O(k^{-1}).
\]
Equation~\eqref{eq:adaagm-curvature-certificate} also gives
$q_k^{\rm AGM}\le q_\infty U_k$.

At a noncoincident index, i.e., $x_k\neq x_{k-1}$, we have $s_kL_k\le q_k^{\rm AGM}\le q_\infty U_k$. 
If $x_k=x_{k-1}\neq x^\star$, then we have $x_{k+1}\neq x_k$. 
In this case, we have $L_{k}=L_{k-1}$ and
\(s_kL_k \le U_ks_{k-1}L_{k-1} \le U_kq_{k-1}^{\rm AGM} \le q_\infty U_kU_{k-1}.
\)
Since $U_{k-1}>1$, we conclude that
\begin{equation*}
    s_kL_k\le q_\infty U_kU_{k-1},
    \qquad k\ge2,
\end{equation*}
unless the full state has already reached the minimizer.

Substituting these bounds into~\eqref{eq:NC-P4} gives
\[
\begin{aligned}
    \chi_k\ge{}&
    b-[1+\gamma b+b^2+bv_{0,k}]q_\infty U_kU_{k-1}-v_{0,k}U_k^2-2U_k(U_k-1).
\end{aligned}
\]
The right-hand side tends to $\delta_{\rm AGM}>0$.
Hence there is a finite index $K_{\rm AGM}\ge2$ such that
\[
    \chi_k\ge\delta_{\rm AGM}/2>0,
    \qquad k\ge K_{\rm AGM}.
\]
Unless the state has already become stationary, Assumption~\ref{ass:NC-sufficient} therefore holds from $K_{\rm AGM}$, with a uniform positive margin in \textup{(E3)}.
The bounds used in the proof of Corollary~\ref{cor:unified-uniform-rate} then imply $\inf_{k\ge K_{\rm AGM}}\widehat\rho_k>0$.
Since $W_k=\Theta(k^2)$, Corollary~\ref{cor:NC-sufficient} yields~\eqref{eq:unified-linear-rate} with $r_\theta=1$ and $K=K_{\rm AGM}$, without assuming stepsize convergence.

\begin{remark}
Fix admissible AdaAGM parameters with $\Omega,\eta\in(0,1)$ and $s_0\ge M_{\rm AGM}/L$.
Our bound admits a geometric multiplier $1-c_{\rm AGM}/\kappa_f$, where $c_{\rm AGM}>0$ depends only on the algorithmic parameters.
This matches the geometric-parameter dependence in
Theorem~4.1 of~\cite{wang2025adaagm}.
\end{remark}

\subsection{Framework-guided algorithm design}
\label{subsec:framework-guided-instances}

We now choose the coefficients and stepsize rule to satisfy the same sufficient conditions directly.
After choosing $\alpha_k$ and $\theta_k$, compute $U_k$ from \textup{(E2)} and require $u_k\le U_k$.
Choose an upper bound $q_k$ for $s_kL_k$ such that
$b_k-(1+\alpha_kb_k+b_k^2)q_k>0$.
The variable-coupling (VC) and constant-coupling (CC) families below follow this construction.

\subsubsection{Admissible families and update rules}
\label{subsec:framework-guided-construction}
Choose $p>3$, $0<\bar\alpha<2$, and
$\theta_k=(k+p)/p$ for $k\ge0$. Write
$b(\bar\alpha):=\min\{1,2-\bar\alpha\}$. Define
\begin{equation}\label{eq:framework-guided-families}
\begin{aligned}
    \alpha_k^{\rm VC}
       &=\bar\alpha\left(\frac{k+1}{k+p}\right)^2,
       &q_k^{\rm VC}&=\frac{2-\alpha_k^{\rm VC}}4,\\
    \alpha_k^{\rm CC}&=\bar\alpha,
       &q_k^{\rm CC}&=
          \nu\frac{b(\bar\alpha)}
                     {1+\bar\alpha b(\bar\alpha)+b(\bar\alpha)^2},
          \qquad 0<\nu<1.
\end{aligned}
\end{equation}
The VC coupling increases to $\bar\alpha$, with
$\alpha_k^{\rm VC}\theta_k^2=\bar\alpha(k+1)^2/p^2$.
For CC, the factor $\nu\in(0,1)$ selects a strict fraction of
$b(\bar\alpha)/(1+\bar\alpha b(\bar\alpha)+b(\bar\alpha)^2)$.
Write the resulting constant as $q^{\rm CC}:=q_k^{\rm CC}$.
Substitution into the definition of $U_k$ in \textup{(E2)} gives
\begin{equation*}
    U_{k+1}^{\rm VC}=\frac{(k+1)(k+p+1)}{(k+2)^2}>1,\qquad
    U_{k+1}^{\rm CC}=\frac{(k+p)^2}{(k+p+1)(k+1)}>1.
\end{equation*}
Both multipliers exceed one for every $k\ge0$ because $p>3$.
For either policy, initialize $x_0=z_0$, obtain
$L_0\in[\mu,L]$ from a probe, and use
\begin{equation}\label{eq:framework-guided-stepsize}
    s_{k+1}=\min\left\{U_{k+1}s_k,
                              \frac{q_{k+1}}{L_{k+1}}\right\},
    \qquad k\ge0,\qquad s_0=\frac1{32L_0},
\end{equation}
together with \eqref{eq:scheme} and \eqref{eq:Lk-bregman-form}.
Here $\theta_0=1$. After the first update, the previous gradient-step and coupling identities hold, so the Lyapunov analysis applies for $k\ge1$.

The update gives $u_k\le U_k$ for $k\ge1$. Since $U_{k+1}>1$ and $L_k\le L$, induction gives
\[
    s_k\ge\frac1L
    \begin{cases}
       \min\{1/32,(2-\bar\alpha)/4\},&\mathrm{VC},\\[1mm]
       \min\{1/32,q^{\rm CC}\},&\mathrm{CC},
    \end{cases}
    \qquad k\ge0.
\]
For every $k\ge0$, $\bar\alpha/p^2\le\alpha_k^{\rm VC}<\bar\alpha<2$ and $\alpha_k^{\rm CC}=\bar\alpha\in(0,2)$. Together with $\theta_k>1$ for $k\ge1$, these bounds verify \textup{(E1)}.
For $k\ge1$, the curvature restriction gives
\begin{equation}\label{eq:framework-guided-curvature-margins}
    b_k-(1+\alpha_kb_k+b_k^2)s_kL_k
    \ge
    \begin{cases}
       \alpha_k^2/4,&\mathrm{VC},\ \alpha_k\le1,\\[1mm]
       (2-\alpha_k)(2\alpha_k-1)/4,&\mathrm{VC},\ \alpha_k>1,\\[1mm]
       (1-\nu)b(\bar\alpha),&\mathrm{CC}.
    \end{cases}
\end{equation}
For each fixed admissible parameter choice, the right-hand sides above have positive lower limits.
Also $\alpha_k\to\bar\alpha$, $v_{0,k}=O(k^{-1})$, and $U_k-1=O(k^{-1})$.
The additional terms subtracted in the definition of $\chi_k$ therefore tend to zero, giving
$\liminf\limits_{k\to\infty}\chi_k>0$.
Moreover,
\[
    W_k^{\rm VC}=\frac{\bar\alpha k(k+1)^2}{p^2(k+p)},\qquad
    W_k^{\rm CC}=\frac{\bar\alpha k(k+p)}{p^2},
\]
Thus both families satisfy $W_k=\Theta(k^2)$ and the growth order used in \eqref{eq:unified-linear-rate} with $r_\theta=1$.

\subsubsection{Selected instances and their guarantee}
\label{subsec:framework-guided-selected}
The tested instances use
\begin{equation}\label{eq:framework-guided-selected-choices}
\begin{aligned}
    \mathrm{VC}:&\quad(p,\bar\alpha)=(20,48/25),\\
    \mathrm{CC}:&\quad(p,\bar\alpha,\nu)=(26,7/5,19/20).
\end{aligned}
\end{equation}
These values were selected by numerical exploration.
Both instances use \eqref{eq:framework-guided-stepsize}; the corollary below proves their eventual linear convergence.

For these selections,
$\alpha_0^{\rm VC}=3/625$ and $\alpha_0^{\rm CC}=7/5$, while for every $k\ge0$,
\[
    q_k^{\rm VC}>\frac1{50},\qquad
    q^{\rm CC}=\frac{57}{220},\qquad
    s_k^{\rm VC}\ge\frac1{50L},\qquad
    s_k^{\rm CC}\ge\frac1{32L}.
\]
For every $k\ge4$, inequality \eqref{eq:unified-delta-margin} holds with
\[
    \delta_{\rm VC}=\frac1{576},\qquad
    \delta_{\rm CC}=\frac3{100}.
\]
Indeed, $\alpha_k^{\rm VC}\ge\alpha_4^{\rm VC}=1/12$ for $k\ge4$.
If $\alpha_k^{\rm VC}\le1$, the first bound in
\eqref{eq:framework-guided-curvature-margins} is at least $1/576$.
On $1\le\alpha_k^{\rm VC}\le48/25$, the second bound is a concave quadratic; its endpoint values are $1/4$ and $71/1250$, both exceeding $1/576$.
For CC, the lower bound is exactly
$(1-\nu)b(\bar\alpha)=(1/20)(3/5)=3/100$.
Lemma~\ref{lem:unified-coefficient-bounds} shows the bounds on coefficient changes.

\begin{corollary}[Selected framework-guided instances]
\label{cor:framework-guided-instances}
Let $f\in\mathcal S^1_{\mu,L}$.
For each selected method in \eqref{eq:framework-guided-selected-choices}, initialized by \eqref{eq:framework-guided-stepsize}, there is a finite index $K\ge4$ such that Assumption~\ref{ass:NC-sufficient} holds for every $k\ge K$.
For each fixed problem and initialization, there is $c>0$ with
$0<c/\kappa_f<1$ such that
\[
    f(x_k)-f_\star,\quad\norm{x_k-x^\star}^2
       =O\!\left(k^{-2}(1-c/\kappa_f)^{k-K}\right).
\]
\end{corollary}

\begin{proof}
The lower bounds on $s_k$ and $\alpha_k$ give \textup{(E1)}, and the stepsize update gives \textup{(E2)} for $k\ge1$.
For $k\ge4$, the constants $\delta_{\rm VC}$ and $\delta_{\rm CC}$ establish \eqref{eq:unified-delta-margin}.
Together with Lemma~\ref{lem:unified-coefficient-bounds}, these verify the hypotheses of Theorem~\ref{thm:unified-tail-verification}.
Since $W_k=\Theta(k^2)$, Corollary~\ref{cor:unified-uniform-rate} gives the stated rate.
\end{proof}

\section{Numerical experiments}
\label{sec:numerics}

We compare seven methods on four randomly generated problem families, measuring both relative errors and the iterations required to reduce them to prescribed tolerances.

All problems have dimension $d=32$.
Using the global bounds $\mu$ and $L$ specified below, we prescribe
\[
    \kappa_f=L/\mu\in\{2,10,10^2,10^3,10^4,10^5,10^6\}.
\]

\begin{itemize}
	\item Ridge regression. Take $n=24$ and draw $A\in\R^{n\times d}$ with independent
	$\mathrm{Unif}[0,1]$ entries and, independently, draw $u$ uniformly from
	the Euclidean unit ball. Set $b=Au$ and define
	\[
	f(x)=\frac{1}{2n}\norm{Ax-b}^2+\frac{\lambda}{2}\norm{x}^2 .
	\]
	We have  $x^\star=(A^\top A+n\lambda I)^{-1}A^\top b$, $L=\norm{A}_2^2/n+\lambda$ and $\mu=\lambda$.
	\item Asymmetric squared loss. Draw $b\sim\mathcal N(0,I_d)$. With $t_\pm=\max(\pm t,0)$, set
	\[
	f(x)=\frac12\sum_{j=1}^d
	\left[\kappa_f(x_j-b_j)_+^2+(x_j-b_j)_-^2\right].
	\]
	We have $x^\star=b$, $f_\star=0$, $\mu=1$ and $L=\kappa_f$.
	\item Binary logistic regression. Take $n=128$ and draw $A\in\R^{n\times d}$ with independent standard
	normal entries. Give each row $a_i^\top$ an independent sign
	$y_i\in\{-1,1\}$ and define
	\[
	f(x)=\frac1n\sum_{i=1}^n\log(1+\exp(-y_i a_i^\top x))
	+\frac{\lambda}{2}\norm{x}^2 .
	\]
	We have $L=\norm{A}_2^2/(4n)+\lambda$ and $\mu=\lambda$.
	\item Log-sum-exp. Draw $128$ affine functions independently, each with slope
	$a_i\sim\mathcal N(0,I_d)$ and independent offset $r_i\sim\mathcal N(0,1)$.
	Define
	\[
	f(x)=\log\sum_{i=1}^{128}\exp(a_i^\top x-r_i)+\frac{\lambda}{2}\norm{x}^2 .
	\]
	We have $L=\max_{i,j}\norm{a_i-a_j}^2/4+\lambda$ and $\mu=\lambda$.
\end{itemize}

For each family we generate five independent draws. On each draw we hold all sampled quantities, including one independent starting point $x_0\sim\mathcal N(0,I_d)$, fixed across every method and conditioning level. 
Changing $\kappa_f$ varies only the parameters in the problem families.
This gives $4\times5\times7=140$ problems.

AdaNAG-G$_{12}$, AdaNAG-G$_{1/2}$, AC-FGM, AdaAGM$_{5.5}$, VC-FG and
CC-FG use the coefficients in Section~\ref{sec:framework-consequences}.
Here AdaAGM$_{5.5}$ denotes the specialization in Corollary~5.5
of Wang and Peypouquet's original preprint.
AC-FGM denotes the smooth unconstrained specialization of
Li and Lan's stepsize policy II~\cite[Corollary~2]{li2026uniformly},
with $a=1/10$ and $\beta^{\rm AC}=1-\sqrt6/3$.
These six methods initialize $L_0$ using the secant
between $x_0$ and $x_0-g_0/\norm{g_0}$. Oracle NAG-SC receives $\mu$ and $L$.
The initial point has index $k=0$, and each accepted update increases $k$
by one.

We retain each method's initialization policy.
AdaNAG-G$_{12}$ and AdaNAG-G$_{1/2}$ use $s_0=21/(802L_0)$ and $s_0=2\sqrt3/(5L_0)$, respectively~\cite[Corollaries~7--8]{suh2025adanag}.
For AdaAGM$_{5.5}$, we take $(m,\gamma,\beta^{\rm AGM},\theta_0)=(1/8,1/2,1,2)$, which gives $s_0=1/(12L_0)$ under the initialization in Section~\ref{sec:framework-consequences}; VC-FG and CC-FG use $s_0=1/(32L_0)$.
AC-FGM follows the authors' first-step search:\footnote{\url{https://github.com/tli432/AC-FGM-Implementation}} start at $\eta^{\rm init}=8/(5L_0)$ and multiply by $2/3$ until $\eta^{\rm init}L_1\le2/5$, recomputing the secant $L_1$ at each trial.
Oracle NAG-SC uses the constant stepsize $1/L$. 
The code to reproduce the results is available at
\url{https://github.com/WilkemZhang/sc-framework-experiments}.

Using $\Delta_k=f(x_k)-f_\star$ from Section~\ref{subsec:setting-notation},
we measure
\begin{equation*}
 \frac{\Delta_k}{\Delta_0},
 \qquad
 \frac{\norm{x_k-x^\star}^2}{\norm{x_0-x^\star}^2}.
\end{equation*}
Both ratios equal one at $k=0$.
For the ridge regression and asymmetric squared loss, we can find the exact solution $x^\star$ and $f_\star$.
For the other two families, we use L-BFGS-B with Newton refinement when needed for binary logistic
and log-sum-exp to approximately find $x^\star$ and $f_\star$. 

\begin{figure}[!htb]\centering
\includegraphics[width=0.97\linewidth]{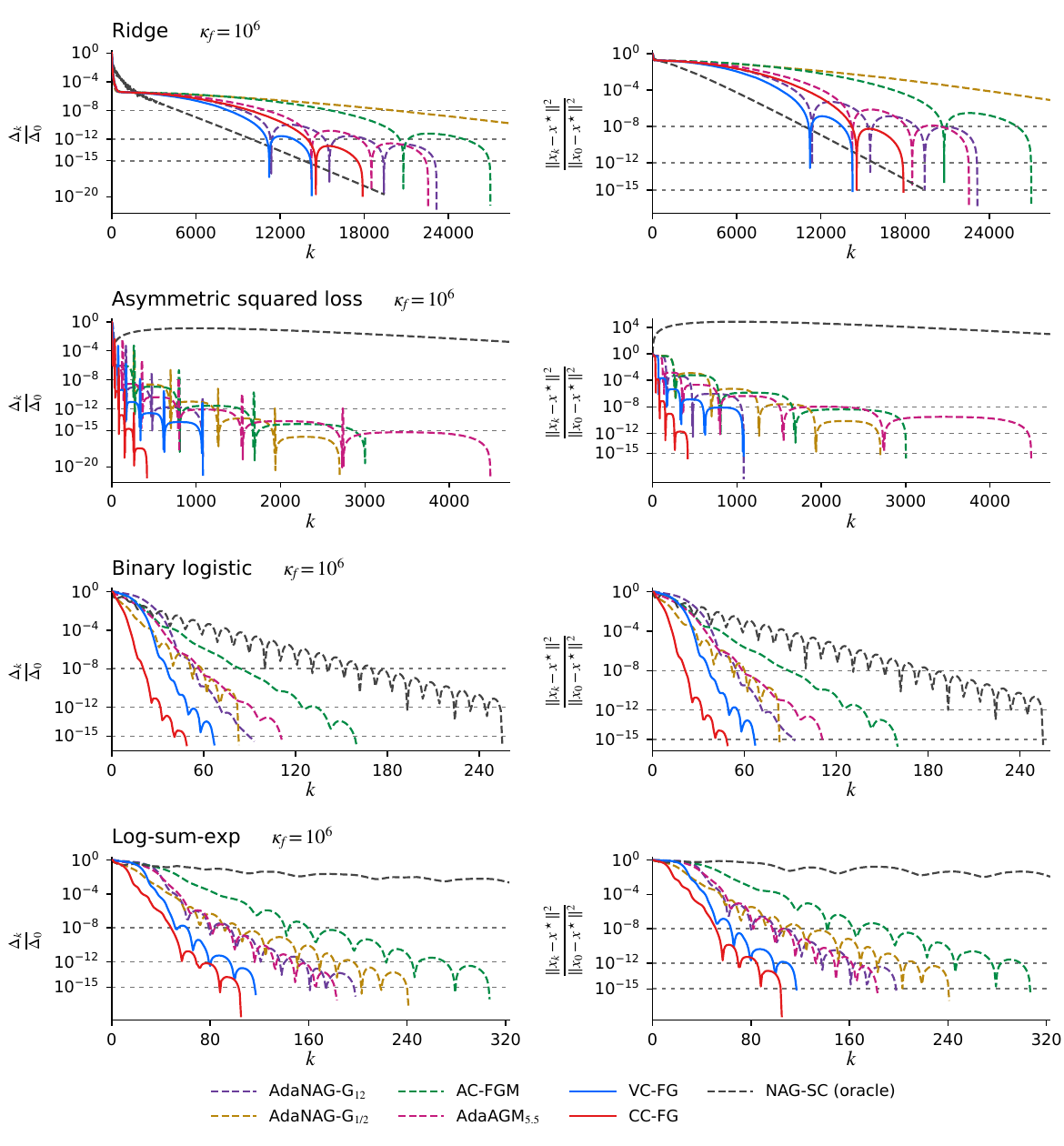}
\caption{Relative objective errors (left) and relative squared-distance errors (right) for the first draw of each problem family at $\kappa_f=10^6$. Rows show ridge regression, asymmetric squared loss, binary logistic regression, and log-sum-exp, respectively. 
Horizontal dashed lines mark $10^{-8}$, $10^{-12}$, and $10^{-15}$.}
\label{fig:num-all-families}
\end{figure}

Figure~\ref{fig:num-all-families} shows the relative objective and squared-distance errors for the first draw of each family at $\kappa_f=10^6$.
On the asymmetric squared loss, binary logistic regression, and log-sum-exp examples, all six adaptive methods reach the joint $10^{-15}$ target before oracle NAG-SC, with CC-FG requiring the fewest iterations. 
VC-FG is fastest on the ridge example, where the separation between methods is less pronounced.
The trajectories also exhibit oscillations and different crossing times for the two errors: a temporary decrease in the objective error need not coincide with a sufficiently small squared-distance error.
This illustrates the importance of requiring both targets at the same iterate.

\begin{figure}[!htb]
\centering
\includegraphics[width=0.95\linewidth]{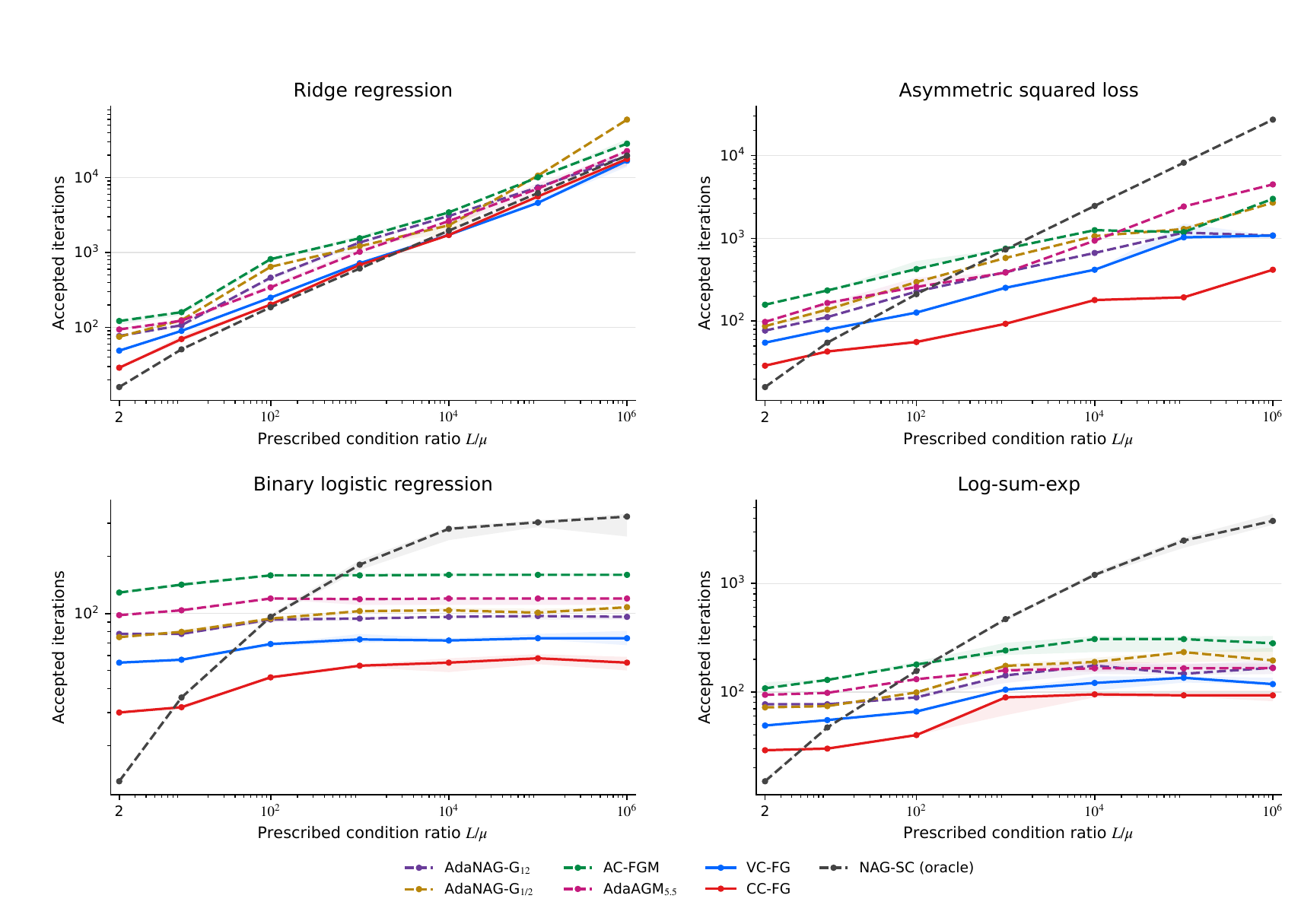}
\caption{Median iteration counts at joint relative tolerance $10^{-15}$.
Lines show the median accepted iterations required for both relative errors to be at most $10^{-15}$ at the same iterate, over five matched
draws. Shaded bands show the interquartile ranges.}
\label{fig:num-median-iterations}
\end{figure}

To complement the trajectories of individual instances, we compare the median iterations required to reach the joint \(10^{-15}\) target over five matched draws at each conditioning level. 
Figure~\ref{fig:num-median-iterations} shows that the adaptive methods perform particularly well on asymmetric squared loss, binary logistic regression, and log-sum-exp, especially at large prescribed condition ratios.
For every tested $\kappa_f\ge10^4$, all six adaptive methods require fewer median iterations than oracle NAG-SC on these three families.
On these three families, CC-FG has the smallest median for every tested $\kappa_f\ge10$.
For binary logistic regression and log-sum-exp, the adaptive iteration counts vary little over the largest conditioning levels, while those of oracle NAG-SC continue to increase.
Ridge regression exhibits a different pattern: iteration counts grow substantially with $\kappa_f$, and the adaptive methods do not provide a uniform advantage over the oracle baseline, although VC-FG and CC-FG remain competitive.

\section{Conclusion}\label{sec:conclusion}

We developed a common Lyapunov analysis for adaptive accelerated gradient methods on smooth, strongly convex objectives.
For stepsizes converging to a positive limit, asymptotic coefficient conditions yield eventual contraction.
Without assuming stepsize convergence, we derive iteration-dependent contraction estimates.
Explicit conditions on local curvature, ratios of consecutive stepsizes, and coefficient changes ensure that the contraction factors eventually remain below a fixed constant smaller than one, giving eventual $R$-linear convergence.

We verified these conditions for specified AdaNAG-G, smooth unconstrained AC-FGM, and AdaAGM policies and used them to construct adaptive families with variable or constant coupling.
The fixed-step FISTA/NAG-type specialization establishes eventual linear convergence for every $0<s\le1/L$.
For $0<s<1/L$, it improves the exponential factor in the stated Li--Shi--Yuan bound for the same smooth iteration.

Experiments on four strongly convex problem families compare complete adaptive policies and illustrate complementary performance of the selected VC-FG and CC-FG methods.


%
%

\bibliographystyle{spmpsci}      
\bibliography{introduction_refs}   

\begin{appendices}

\section{Coefficient growth and curvature bounds}
\label{app:section6-coefficient-bounds}

We derive the coefficient bounds and constants used in
Theorem~\ref{thm:unified-tail-verification}, together with the growth
estimate for $W_k$ used in Corollary~\ref{cor:unified-uniform-rate}.
We then prove the bounds on $q_k^G$ and $U_k$ used to verify the two
AdaNAG-G policies in Section~\ref{sec:framework-consequences}.
\subsection{Common coefficient policies}
\label{subsubsec:unified-coefficient-sequences}

For \(k\ge2\), the represented coefficient policies have the form
\begin{equation}
\label{eq:unified-coefficient-representation}
\begin{aligned}
 \theta_{k+1}-\theta_k
 &=c_{\Delta\theta,k}(k+1)^{r_\theta-1},\\
 \alpha_k
 &=c_\alpha
 \left(\frac{\theta_{k+1}-d_{\alpha,1}}
 {\theta_k-d_{\alpha,2}}\right)^{r_\alpha}.
\end{aligned}
\end{equation}
The quantities $r_\theta,r_\alpha,c_\alpha,d_{\alpha,1},d_{\alpha,2}$
are fixed, while $c_{\Delta\theta,k}$ may vary with $k$.
The exponent $r_\theta$ specifies the power in the increment formula;
$r_\alpha=0$ gives constant coupling.
Table~\ref{tab:unified-coefficient-representation} lists the
parameters for the methods in Section~\ref{sec:framework-consequences}.
For AC-FGM, $a\in(0,1]$ ensures a positive lower bound on the increments.

\begin{table}[ht]
\centering
\small
\caption{Parameters in the common coefficient representation
\eqref{eq:unified-coefficient-representation} and bounds
\(c_{\Delta\theta,k}\in[\underline c_\theta,\overline c_\theta]\).}
\label{tab:unified-coefficient-representation}
\renewcommand{\arraystretch}{1.20}
\setlength{\tabcolsep}{4pt}
\begin{tabularx}{\textwidth}{>{\raggedright\arraybackslash}X c c c c c}
\toprule
Method & \(r_\theta\) & \([\underline c_\theta,\overline c_\theta]\)
& \(c_\alpha\) & \((d_{\alpha,1},d_{\alpha,2})\) & \(r_\alpha\)\\
\midrule
AdaNAG-G\(_{12}\) & \(1\) & \([1/12,1/12]\)
& \(1/2\) & \((1,0)\) & \(2\)\\
AdaNAG-G\(_{1/2}\) & \(1/2\) & \([2/(\sqrt5+2),1]\)
& \(1/2\) & \((0,0)\) & \(0\)\\
Adaptive-\(\tau\) AC-FGM & \(1\) & \([a/2,1/2]\)
& \(\beta^{\rm AC}/(1-\beta^{\rm AC})\) & \((1,0)\) & \(1\)\\
AdaAGM & \(1\) & \([m/2,m]\)
& \(\gamma\) & \((0,0)\) & \(0\)\\
VC family & \(1\) & \([1/p,1/p]\)
& \(\bar\alpha\) & \((1,0)\) & \(2\)\\
CC family & \(1\) & \([1/p,1/p]\)
& \(\bar\alpha\) & \((0,0)\) & \(0\)\\
\bottomrule
\end{tabularx}
\end{table}

The VC and CC rows describe the families in
\eqref{eq:framework-guided-families}; the selected VC-FG and CC-FG
parameters are given in \eqref{eq:framework-guided-selected-choices}.
In the CC family, $\nu$ enters the stepsize rule through the bound
$q^{\rm CC}$. All parameter ranges are as stated in
Section~\ref{sec:framework-consequences}.

Assume
\begin{equation}\label{eq:unified-coefficient-admissibility}
\begin{gathered}
 0<\underline c_\theta\le c_{\Delta\theta,k}
 \le\overline c_\theta\le1,\qquad
 0<r_\theta\le1,\qquad \theta_2>1,\\
 c_\alpha>0,\qquad r_\alpha\ge0,\qquad
 0\le d_{\alpha,1},d_{\alpha,2}\le1.
\end{gathered}
\end{equation}
Section~\ref{sec:framework-consequences} separately verifies the range
$\alpha_k<2$, the stepsize and coefficient conditions
\textup{(E1)}--\textup{(E2)},
and the positive lower bound in \eqref{eq:unified-delta-margin}.

\subsection{Coefficient growth and variation}
\label{subsubsec:unified-growth-orders}

Write
\[
    \Delta_k^\theta:=\theta_{k+1}-\theta_k,\qquad
    J_k:=\frac{\theta_{k+1}-d_{\alpha,1}}
                 {\theta_k-d_{\alpha,2}},
\]
so that $\alpha_k=c_\alpha J_k^{r_\alpha}$.
First define constants for the bounds on $\theta_k$, $J_k$, and
$\alpha_k$:
\begin{equation}
\label{eq:unified-growth-constants}
\begin{aligned}
 D_\theta&:=\frac{r_\theta}{\underline c_\theta
 [1-(3/4)^{r_\theta}]},\\
 J_-&:=\frac{3^{r_\theta}}{D_\theta+3^{r_\theta}},
 &J_+&:=1+\frac{(1+\overline c_\theta)D_\theta}{3^{r_\theta}},\\
 \alpha_-&:=c_\alpha J_-^{r_\alpha},
 &\alpha_+&:=c_\alpha J_+^{r_\alpha}.
\end{aligned}
\end{equation}
Next, define constants for the differences between consecutive
ratios $J_k$ and coupling coefficients $\alpha_k$:
\begin{equation}
\label{eq:unified-increment-constants}
\begin{aligned}
 B_R:={}&D_\theta[\overline c_\theta-\underline c_\theta
 +\overline c_\theta(1-r_\theta)]\\
 &+\overline c_\theta
 (\overline c_\theta+|d_{\alpha,1}-d_{\alpha,2}|)D_\theta^2,\\
 A_\alpha:={}&\frac{r_\alpha B_R}{J_-^{r_\alpha}}
 \max\{J_-^{r_\alpha-1},J_+^{r_\alpha-1}\}.
\end{aligned}
\end{equation}
For constant coupling, \(A_\alpha=0\).
Finally, combine these constants to obtain the bounds on
$v_{0,k}$ and $U_k$:
\begin{equation}
\label{eq:unified-E4-coefficient}
\begin{aligned}
 C_0:={}&(1+A_\alpha)(1+\overline c_\theta D_\theta)^2-1
 +\alpha_-A_\alpha+\alpha_+D_\theta,\\
 C_U:={}&(1+2^{r_\theta}A_\alpha)(1+D_\theta)-1,
 &\overline U:={}&1+C_U.
\end{aligned}
\end{equation}

\begin{lemma}[Common coefficient bounds]
\label{lem:unified-coefficient-bounds}
Under \eqref{eq:unified-coefficient-representation}--
\eqref{eq:unified-coefficient-admissibility}, for every \(k\ge3\),
\begin{equation}
\label{eq:unified-coefficient-growth-bound}
 \frac1{\theta_k-1}\le\frac{D_\theta}{k^{r_\theta}},
 \qquad
 \frac{\theta_{k+1}-\theta_k}{\theta_k-1}
 \le\frac{\overline c_\theta D_\theta}{k},
\end{equation}
and
\begin{equation}\label{eq:unified-coefficient-orders}
 \alpha_-\le\alpha_k\le\alpha_+,\qquad
 \theta_k-1=\Theta(k^{r_\theta}),\qquad
 W_k=\Theta(k^{2r_\theta}).
\end{equation}
For every \(k\ge4\),
\begin{equation}
\label{eq:unified-tail-constants}
 v_{0,k}\le\frac{C_0}{k^{r_\theta}},\qquad
 0<U_k\le1+\frac{C_U}{k^{r_\theta}}\le\overline U,
 \qquad
 [U_k-1]_+\le\frac{C_U}{k^{r_\theta}}.
\end{equation}
The bounds on \(U_k\) do not require \(U_k\ge1\).
\end{lemma}

\begin{proof}
\emph{Growth of $\theta_k$.}
For $k\ge3$, summing the lower increment bound gives
\[
\begin{aligned}
    \theta_k-1
    &\ge\underline c_\theta
          \sum_{j=2}^{k-1}(j+1)^{r_\theta-1}\\
    &\ge\underline c_\theta\int_3^{k+1}x^{r_\theta-1}\,dx
     =\frac{\underline c_\theta}{r_\theta}
       [(k+1)^{r_\theta}-3^{r_\theta}]\ge\frac{k^{r_\theta}}{D_\theta}.
\end{aligned}
\]
The last inequality uses
$(k+1)^{r_\theta}-3^{r_\theta}
 \ge k^{r_\theta}[1-(3/4)^{r_\theta}]$ and the definition of $D_\theta$.
Summing the upper increment bound similarly gives
$\theta_k-1=O(k^{r_\theta})$.
Since \((k+1)^{r_\theta-1}\le k^{r_\theta-1}\), these estimates also
give \eqref{eq:unified-coefficient-growth-bound}.

\emph{Bounds on the coupling coefficients.}
The restrictions on $d_{\alpha,1},d_{\alpha,2}$ give
\begin{equation}\label{eq:unified-R-bounds}
 \frac{\theta_k-1}{\theta_k}\le J_k
 \le\frac{\theta_k+\Delta_k^\theta}{\theta_k-1},
 \qquad J_-\le J_k\le J_+,\qquad k\ge3.
\end{equation}
Indeed, the preceding bounds on $\theta_k$ imply
\[
    \frac1{1+D_\theta/k^{r_\theta}}
    \le J_k
    \le1+\frac{(1+\overline c_\theta)D_\theta}{k^{r_\theta}}.
\]
Using $k\ge3$ gives the constants $J_-$ and $J_+$ above.
Thus $\alpha_k=c_\alpha J_k^{r_\alpha}$ lies between $\alpha_-$
and $\alpha_+$. Together with $\theta_k-1=\Theta(k^{r_\theta})$,
these bounds give
$W_k=\alpha_k\theta_k(\theta_k-1)=\Theta(k^{2r_\theta})$,
proving \eqref{eq:unified-coefficient-orders}.

\emph{Differences between consecutive coefficients.}
Subtracting the ratios gives
\[
 J_{k+1}-J_k
 =\frac{\Delta_{k+1}^\theta-\Delta_k^\theta}
 {\theta_{k+1}-d_{\alpha,2}}
 -\frac{(\Delta_k^\theta-d_{\alpha,1}+d_{\alpha,2})
 \Delta_k^\theta}
 {(\theta_k-d_{\alpha,2})(\theta_{k+1}-d_{\alpha,2})}.
\]
The increment policy gives
\[
 |\Delta_{k+1}^\theta-\Delta_k^\theta|
 \le(\overline c_\theta-\underline c_\theta)(k+1)^{r_\theta-1}
 +\overline c_\theta(1-r_\theta)(k+1)^{r_\theta-2}.
\]
Use \(\theta_j-d_{\alpha,2}\ge j^{r_\theta}/D_\theta\),
\(\Delta_k^\theta\le\overline c_\theta\), and
\(|\Delta_k^\theta-d_{\alpha,1}+d_{\alpha,2}|
 \le\overline c_\theta+|d_{\alpha,1}-d_{\alpha,2}|\).
The two terms in the ratio difference therefore satisfy
\[
\begin{aligned}
    |J_{k+1}-J_k|
    \le{}&D_\theta\left[
       \frac{\overline c_\theta-\underline c_\theta}{k+1}
       +\frac{\overline c_\theta(1-r_\theta)}{(k+1)^2}
    \right]\\
    &+\frac{\overline c_\theta
       (\overline c_\theta+|d_{\alpha,1}-d_{\alpha,2}|)D_\theta^2}
       {k^{r_\theta}(k+1)^{r_\theta}}
    \le\frac{B_R}{k^{r_\theta}}.
\end{aligned}
\]
The last step uses $0<r_\theta\le1$.
If $r_\alpha=0$, the coupling is constant and $A_\alpha=0$.
If $r_\alpha>0$, apply the mean-value theorem to $J^{r_\alpha}$
on $[J_-,J_+]$. Since
\[
    \alpha_-A_\alpha
    =c_\alpha r_\alpha B_R
       \max\{J_-^{r_\alpha-1},J_+^{r_\alpha-1}\},
\]
both cases give
\begin{equation}\label{eq:unified-alpha-increment-bound}
 |\alpha_{k+1}-\alpha_k|
 \le\frac{\alpha_-A_\alpha}{k^{r_\theta}},\qquad
 \frac{|\alpha_{k+1}-\alpha_k|}{\alpha_k}
 \le\frac{A_\alpha}{k^{r_\theta}},\qquad k\ge3.
\end{equation}

\emph{Bounds on $v_{0,k}$ and $U_k$.}
Use the identities
\begin{equation}\label{eq:unified-proof-quantities}
\begin{aligned}
 \frac{W_{k+1}}{W_k}
 &=\left(1+\frac{\alpha_{k+1}-\alpha_k}{\alpha_k}\right)
 \left(1+\frac{\Delta_k^\theta}{\theta_k}\right)
 \left(1+\frac{\Delta_k^\theta}{\theta_k-1}\right),\\
 U_k&=\frac{\alpha_{k-1}}{\alpha_k}
 \left(\frac{\theta_{k-1}}{\theta_k}\right)^2
 \left(1+\frac1{\theta_k-1}\right).
\end{aligned}
\end{equation}
Expand the first product and bound each term by its absolute value.
Using \eqref{eq:unified-alpha-increment-bound},
\eqref{eq:unified-coefficient-growth-bound}, and
$k^{-1}\le k^{-r_\theta}$ gives
\[
\begin{aligned}
    \left|\frac{W_{k+1}}{W_k}-1\right|
    &\le\left(1+\frac{A_\alpha}{k^{r_\theta}}\right)
          \left(1+\frac{\overline c_\theta D_\theta}{k}\right)^2-1\\
    &\le\frac{(1+A_\alpha)(1+\overline c_\theta D_\theta)^2-1}
                   {k^{r_\theta}}.
\end{aligned}
\]
Together with \eqref{eq:unified-alpha-increment-bound} and
\(\alpha_k/(\theta_k-1)\le\alpha_+D_\theta/k^{r_\theta}\), this
proves the bound on \(v_{0,k}\).
For the second identity, the estimate on
$\alpha_k-\alpha_{k-1}$ uses
\eqref{eq:unified-alpha-increment-bound} at index $k-1$.
We therefore take $k\ge4$, so that $k-1\ge3$. Then
\[
 \frac{\alpha_{k-1}}{\alpha_k}
 \le1+\frac{\alpha_-A_\alpha}{\alpha_k(k-1)^{r_\theta}}
 \le1+\frac{2^{r_\theta}A_\alpha}{k^{r_\theta}}.
\]
Since $\theta_{k-1}/\theta_k<1$, substituting this bound into
the identity for $U_k$ gives
\[
 U_k\le
 \left(1+\frac{2^{r_\theta}A_\alpha}{k^{r_\theta}}\right)
 \left(1+\frac{D_\theta}{k^{r_\theta}}\right)
 \le1+\frac{C_U}{k^{r_\theta}}.
\]
All factors defining $U_k$ are positive, so $U_k>0$.
Taking the positive part of $U_k-1$ proves the remaining bound in
\eqref{eq:unified-tail-constants}.
\end{proof}

\subsection{Bounds for the two AdaNAG-G policies}
\label{app:section6-native-g-checks}

For the policies in \eqref{eq:G-concrete-policies}, we prove
\[
    \frac1{10}\le q_k^G<\frac13,\qquad U_k>1,\qquad k\ge4.
\]
Since $A_{k-1}^G=W_k$, definition
\eqref{eq:G-general-curvature-certificate} gives
\[
    (q_k^G)^{-1}
    =\frac{W_k}{B_{k-1}^G}
     +\frac{B_k^G+\alpha_k^2\theta_k^2}{W_k}.
\]
The closed forms below use the coefficient formula for $B_k^G$ only
when $k\ge1$; $B_0^G$ retains its separately specified initialization
value. The corresponding closed forms for $q_k^G$ apply for $k\ge2$.

\paragraph{The linear-growth policy.}
For AdaNAG-G$_{12}$, direct substitution gives for $k\ge1$
\[
    W_k=\frac{(k+3)^2(k+2)}{288(k+14)},\qquad
    B_k^G=\alpha_k^2\theta_k^2
           =\frac{(k+3)^4}{576(k+14)^2}.
\]
Consequently, for $k\ge2$,
\[
    q_k^G=\frac{(k+2)^3(k+14)}
                  {(k+3)^2(3k^2+56k+342)}.
\]
The denominator minus three times the numerator is
$14k^3+417k^2+2028k+2742>0$, proving $q_k^G<1/3$.
For the lower bound, ten times the numerator minus the denominator is
\[
    P(k):=7k^4+126k^3+255k^2-796k-1958.
\]
Writing $k=4+j$ gives
\[
    P(4+j)=7j^4+238j^3+2439j^2+9084j+8794>0,
    \qquad j\ge0.
\]
Hence $q_k^G>1/10$ for $k\ge4$.
The upper bound $U_k$ on the stepsize ratio is
\[
    U_k=\frac{(k+2)(k+14)}{(k+3)^2}
        =1+\frac{10k+19}{(k+3)^2}>1.
\]

\paragraph{The square-root-growth policy.}
For AdaNAG-G$_{1/2}$, put $x:=\sqrt{k+3}$.
For $k\ge1$, the coefficient formulas give
\[
    W_k=x(2x-1),\qquad \alpha_k^2\theta_k^2=x^2,\qquad
    B_k^G=B(x):=\frac{x^2(2x^2-4x+3)}{2(x^2-1)}.
\]
For every $t\ge2$,
\[
    \frac{t^2}{2}\le B(t)\le t^2-t\le t^2.
\]
The lower bound follows from
\[
    B(t)-\frac{t^2}{2}
    =\frac{t^2(t-2)^2}{2(t^2-1)}\ge0.
\]
The upper bound follows from
\[
    B(t)-(t^2-t)
       =-\frac{t(2t-1)(t-2)}{2(t^2-1)}\le0.
\]
For $k\ge4$, both $x$ and $\sqrt{x^2-1}$ exceed $2$.
Applying the upper bound on $B$ at $\sqrt{x^2-1}$ gives
\[
    B_{k-1}^G
    \le x^2-1-\sqrt{x^2-1}
    \le x^2-x,
\]
where the last inequality uses $\sqrt{x^2-1}\ge x-1$.
Direct simplification also gives
\[
    \frac{B_k^G+x^2}{W_k}
       =\frac{x(2x-1)}{2(x^2-1)}.
\]
Therefore
\[
    (q_k^G)^{-1}
       \ge\frac{2x-1}{x-1}
          +\frac{x(2x-1)}{2(x^2-1)}
       =3+\frac{x+4}{2(x^2-1)}>3.
\]
For the reverse bound, use $B_{k-1}^G\ge(x^2-1)/2$,
$B_k^G\le x^2$, and $W_k\ge3x^2/2$ to obtain
\[
    (q_k^G)^{-1}
       \le\frac{4x^2}{x^2-1}+\frac43
       \le6<10.
\]
The last inequality uses $x^2=k+3\ge7$ and proves the shared
lower bound $q_k^G\ge1/10$.
The upper bound $U_k$ on the stepsize ratio is
\[
    U_k=\frac{2(k+2)}{2(k+3)-\sqrt{k+3}}>1,
    \qquad k\ge4.
\]

These estimates give the stepsize lower bound and the positive lower
bound in \eqref{eq:unified-delta-margin} used in
Subsubsection~\ref{subsec:G-consequences}.
Taking $k\ge4$ ensures that every coefficient sequence index entering
$q_k^G$ is at least $2$.

\end{appendices}

\end{document}